\documentclass[11pt,leqno]{article}

\usepackage{amsmath,amsfonts,amscd,amssymb,theorem}

\long\def\comment#1\endcomment{}

\makeatletter
\begingroup
\gdef\th@dotted{\normalfont\itshape
  \def\@begintheorem##1##2{%
        \item[\hskip\labelsep \theorem@headerfont ##1\ ##2.]}%
\def\@opargbegintheorem##1##2##3{%
   \item[\hskip\labelsep \theorem@headerfont ##1\ ##2\ (##3).]}}
\endgroup
\makeatother

\theoremstyle{dotted}

\newtheorem{theorem}{Theorem}[section]
\newtheorem{lemma}[theorem]{Lemma}
\newtheorem{conj}[theorem]{Conjecture}
\newtheorem{prop}[theorem]{Proposition}
\newtheorem{corr}[theorem]{Corollary}

\makeatletter
\begingroup
\gdef\th@upshape{\normalfont
  \def\@begintheorem##1##2{%
        \item[\hskip\labelsep \theorem@headerfont ##1\ ##2.]}%
\def\@opargbegintheorem##1##2##3{%
   \item[\hskip\labelsep \theorem@headerfont ##1\ ##2\ (##3).]}}
\endgroup
\makeatother

\theoremstyle{upshape}

\newtheorem{defn}[theorem]{Definition}
\newtheorem{remark}[theorem]{Remark}

\makeatletter
\renewcommand{\subsection}{\@startsection{subsection}{2}{0pt}{-3ex
plus -1ex minus -0.2ex}{-2mm plus -0pt minus
-2pt}{\normalfont\bfseries}} \makeatother

\makeatletter
\@addtoreset{equation}{section}
\makeatother

\newcommand{\cntrct}                
{\hspace{2pt}\raisebox{1pt}{\text{$\lrcorner$}}\hspace{2pt}}

\newcommand{\proof}[1][Proof.]{\smallskip\noindent{\em #1}}
\def\endproof{\hfill\ensuremath{\square}\par\medskip}

\def\eqref#1{\thetag{\ref{#1}}}

\let\latexref=\ref
\def\ref#1{{\normalfont{\latexref{#1}}}}

\newcommand{\wt}{\widetilde}
\newcommand{\wh}{\widehat}

\newcommand{\idot}{{\:\raisebox{1pt}{\text{\circle*{1.5}}}}}
\newcommand{\hdot}{{\:\raisebox{3pt}{\text{\circle*{1.5}}}}}
\newcommand{\tw}{ {(1)} }
\newcommand{\Fr}{{\sf Fr}}
\newcommand{\ad}{\operatorname{\sf ad}}
\newcommand{\cchar}{\operatorname{\sf char}}
\newcommand{\Pic}{\operatorname{\sf Pic}}
\newcommand{\Br}{\operatorname{\sf Br}}
\newcommand{\ZZ}{{\cal Z}}
\newcommand{\Spec}{\operatorname{Spec}}
\newcommand{\sspec}{\operatorname{{\cal S}{\it pec}}}

\newcommand{\Coh}{\operatorname{Coh}}
\newcommand{\fmod}{{\text{\rm -mod}^{\text{{\tt\tiny fg}}}}}

\newcommand{\Z}{{\mathbb Z}}

\newcommand{\HH}{{\sf H}}

\newcommand{\M}{{\cal M}}

\newcommand{\T}{{\cal T}}

\newcommand{\G}{{\sf G}}
\newcommand{\GG}{{\mathcal G}}

\newcommand{\calo}{{\cal O}}

\newcommand{\m}{{\mathfrak m}}

\newcommand{\hh}{{\cal H}}
\newcommand{\gm}{{\mathbb{G}_m}}
\newcommand{\ga}{{\mathbb{G}_a}}

\newcommand{\eps}{\varepsilon}
\renewcommand{\phi}{\varphi}

\renewcommand{\dim}{\operatorname{\sf dim}}
\newcommand{\gr}{\operatorname{\sf gr}}
\newcommand{\id}{\operatorname{\sf id}}
\newcommand{\rk}{\operatorname{\sf rk}}

\newcommand{\Hom}{\operatorname{Hom}}
\newcommand{\Ext}{\operatorname{Ext}}

\newcommand{\Aut}{\operatorname{{\sf Aut}}}

\newcommand{\RR}{{\bf R}}

\newcommand{\loc}{\operatorname{{\sf Loc}}}

\title{Fedosov quantization in positive characteristic}

\author{R. Bezrukavnikov\thanks{Partially supported by NSF grant
DMS-0071967} and D. Kaledin\thanks{Partially supported by CRDF grant
RM1-2694-MO05.}}

\begin{document}

\maketitle

\tableofcontents

\section*{Introduction}

The so-called {\em deformation quantization} of symplectic manifolds
originally appeared in the $C^\infty$ symplectic geometry; the
crucial breakthrough has been made in the early 80-ies independently
by B. Fedosov and M. De Wilde-P. Lecomte. The input of a deformation
quantization problem is a symplectic (or, more generally, a Poisson)
manifold $M$; the output is a non-commutative one-parameter
deformation of the algebra of function on $M$. Both Fedosov and De
Wilde-Lecomte provided general procedures which solve the
deformation quantization problem for $C^\infty$ manifolds. Recently,
motivated in part by M. Kontsevich \cite{Kn}, there was much
interest in generalizing the deformation quantization procedures to
the case of algebraic manifolds equipped with an algebraic symplectic
form (see e.g. \cite{BK}, \cite{Y}, and an earlier paper \cite{NT}
for the holomorphic situation). In particular, in \cite{BK} it has
been shown that under some mild assumptions on the manifold, the
Fedosov quantization can be made to work in the algebraic setting.

The present paper is a continuation of \cite{BK} (to which we refer
the reader for a more complete bibliography and historical
discussion). Namely, one of the most important assumptions in
\cite{BK} (as well as in other papers on the subject) was that the
field of definition for all algebraic manifolds has characteristic
$0$. In this paper, we study what happens in the case of positive
characteristic. 

The most obvious new feature of the theory in positive
characteristic is the presence of a large Poisson center in the
sheaf $\calo_X$ of functions on a Poisson manifold $X$: since for
every local functions $f,g \in \calo_X$ and any Poisson bracket
$\{-,-\}$ on $\calo_X$ we have $\{f^p,g\}=0$, the image $\calo_X^p
\subset \calo_X$ of the Frobenius map lies in the center of any
Poisson structure. This phenomenon, already observed in \cite{BMR},
allows for interesting applications (see e.g. \cite{BK2}) but makes
the quantization procedures more involved. In this paper, we were
not able to prove any meaningful results for general quantizations
in positive characteristic, and we had to restrict our attention to
a special class of them: the so-called {\em Frobenius-constant}
quantization. Roughly speaking -- the precise definition is
Definition~\ref{fr.const.defn} below -- a quantization is
Frobenius-constant if the Poisson center $\calo_X^p \subset \calo_X$
stays central in the quantized algebra $\calo_h$. 

\medskip

For quantization of this type, we were able to achieve, under mild
assumptions on the manifold $X$, a reasonably complete
classification theorem. In particular, Frobenius-constant
quantizations do exist.

Moreover, as a prerequisite to the study of quantizations, we
investigate to some extent symplectic differential geometry in
positive characteristic. In particular, we introduce a notion of a
{\em restricted Poisson algebra} -- a Poisson version of the
standard notion of a restricted Lie algebra. We also prove a version
of the Darboux Theorem for these algebras
(Proposition~\ref{darboux}; in our setting the usual Darboux Theorem
is false). These results seems to be new, and they might be of
independent interest: even a reader who is not interested in
quantizations at all might find symplectic geometry in positive
characteristic worth her (or his) attention.

\medskip

Our main tool is the (version of the) technique of the so-called
{\em formal geometry}; we have used the same technique in
\cite{BK}. In positive characteristic this turns out to be simpler.
In fact, we could have avoided this entirely by working
systematically in the flat topology on the symplectic manifold $X$;
the use of formal geometry in this context amounts to a choice of a
concrete flat covering of $X$ to trivialize things. We decided to
make this choice because the picture becomes more explicit and more
easily accessible to people who have little experience with flat
topology. Moreover, the use of formal geometry emphasizes the
differential-geometric nature of the subject and the parallels with
the characteristic $0$ case.

Our restriction to Frobenius-constant quantizations is particularly
unfortunate because -- as opposed to the main result in \cite{BK} --
it effectively excludes from consideration compact algebraic
varieties, such as abelian varieties or surfaces of the type
K3. There undoubtedly exists a more general theory which would apply
to these cases as well. Another thing which we leave out entirely in
this paper is questions of mixed characteristic: a complete theory
should probably incorporate both the quantizations of a
characteristic $p$ manifold and its liftings to characteristic
$0$. All this should be the subject of future research.

\medskip

The paper is organized as follows. The first section contains the
necessary definitions and the statements of our main
results. Unfortunately, many of this material is not standard. We
have tried to keep the technical details to absolute minimum, but
Section 1 still takes about half of the paper. Since this was
unavoidable anyway, we have also incorporated some of the shorter
proofs into this section; eventually, Section 1 became
self-contained except for four clearly marked main results -- two
deal with the symplectic side of the story, two other with
quantizations. Those two that deal with things symplectic and Poisson
are proved in Section 2. Then in Section 3, we restrict our
attention to local study of symplectic manifolds (this includes our
Darboux-type theorem for restricted Poisson algebras). Finally, in
Section 4 we use the techniques of formal geometry to globalize
things and to prove the remaining two results which deal with
quantizations.

\subsection*{Acknowledgments.} The paper is a part of the ongoing
project for the systematic use of quantizations in positive
characteristic to obtain information about algebraic symplectic
varieties; the project has been going on already for several years,
and during this time we have benefited from the discussions and help
from many people. Our original encouragement came from the late
A.N. Tyurin, who is very much missed. We would like to thank
A. Beilinson, A. Bondal, V. Drinfeld, T. Ekedal, M. Emerton,
P. Etingof, B. Feigin, M. Finkelberg, V. Ginzburg, D. Huybrechts,
D. Kazhdan, M. Kontsevich, A. Kuznetsov, M. Lehn, Ch. Sorger, D. Van
Straten, M. Verbitsky and V. Vologodsky for much help and much
valuable input. During the preparation of the paper, we had the
opportunity to work at several first-rate mathematical institutions
whose hospitality is gratefully acknowledged; in particular, we
would like to thank Northwestern University, the University of
Chicago, Universit\'e Paris 7, Universit\'e de Nantes, Universit\"at
Gutenberg at Mainz, Hebrew University of Jerusalem and the
Mittag-Laefller Institute in Stockholm.

\section{Statements and definitions.}\label{st.sec}

In this section we give the relevant definitions, and state our main
results. The reader is encouraged to compare everything to
\cite[1.3]{BK}.

\subsection{Definitions.} Throughout the paper, we will say that a
subalgebra $A \subset B$ in an algebra $B$ is {\em central} if it
lies in the center of the algebra $B$, and we will say that an
algebra map $A \to B$ is {\em central} if its image is central in
$B$.

Fix once and for all a field $k$ of characteristic $p > 2$. For any
vector space, algebra, scheme etc. $M$ over $k$ denote by $M^\tw = M
\otimes_k k$ its pullback with respect to the Frobenius map $\Fr:k
\to k$. For a commutative algebra $A$ over $k$, denote by $A^p
\subset A$ the subalgebra generated by $p$-th powers of all elements
$a \in A$. By definition, the Frobenius map $\Fr:A \to A$ factors
through the embedding $A \subset A^\tw$ by means of the relative
Frobenius map $\Fr_k:A^\tw \to A$; the subalgebra $A^p \subset A$ is
the image of the map $\Fr_k$. If the algebra $A$ is reduced, the
relative Frobenius map is injective; in this case $A^p \cong
A^\tw$. As usual, by a derivation $D:A \to M$ with values in an
$A$-module $M$ we will understand a $k$-linear map which satisfies
the Leibnitz rule. By a Frobenius-derivation $D:A \to M$ with values
in an $A^p$-module $M$ we will understand a map which is derivation
with respect to the induced $A^\tw$-module structure on $M$ -- in
other words, we have $D(1)=0$ and
$$
D(ab) = a^pD(b)+b^pD(a).
$$
The following definition is exactly the same as in the case $\cchar
k = 0$.

\begin{defn} By a {\em quantization} $A_h$ of a commutative
$k$-algebra $A$ we will understand an associative flat $k[[h]]$-algebra
$A_h$, complete with respect to the $h$-adic filtration and equipped
with an isomorphism $A_h/h \cong A$.

By a {\em quantization} $\calo_h$ of a scheme $X$ over $k$ we will
understand a Zariski sheaf $\calo_h$ of flat $k[[h]]$-algebras on
$X$, complete with respect to the $h$-adic filtration and equipped
with an isomorphism $\calo_h/h \cong \calo_X$.
\end{defn}

As explained in \cite[1.3]{BK}, quantizations of a commutative
$k$-algebra $A$ are in natural one-to-one correspondence with
quantizations of the affine scheme $X = \Spec A$.  Any quantization
$\calo_h$ of a scheme $X/k$ induces a Poisson bracket in $\calo_X$;
one usually keeps this structure in mind from the start and says
that $\calo_h$ is a quantization of the Poisson scheme $X$. We note
that by definition, for any section $f \in \calo_h$ which lies in the
center $\ZZ_h \subset \calo_h$ of the algebra sheaf $\calo_h$, its
reduction $\overline{f} \in \calo_X$ modulo $h$ lies in the Poisson
center $\ZZ \subset \calo_X$.

\begin{defn}\label{cent.defn}
A quantization $\calo_h$ of a scheme $X/k$ is called {\em central}
if the natural map $\ZZ_h \to \ZZ$ from the center $\ZZ_h \subset
\calo_h$ of the associative algebra $\calo_h$ to the Poisson center
$\ZZ \subset \calo_X$ is surjective.
\end{defn}

We note that the center $\ZZ_h \subset \calo_h$ of a quantization
$\calo_h$ is automatically flat over $k[[h]]$ (this is equivalent to
having no $h$-torsion). Thus $\ZZ_h$ is a commutative one-parameter
deformation of some subalgebra $\ZZ_h/h \subset \ZZ$. The quantization
is central if $\ZZ_h/h$ is the whole Poisson center $\ZZ \subset \calo_X$.

In \cite{BK2} we have used a notion of a {\em Frobenius-constant
quantization}; roughly speaking, a quantization was called
Frobenius-constant if it is central, and moreover, we have an
isomorphism $\ZZ_h \cong \calo_X^p[[h]]$. It turns out that for
general theory, we need a more precise notion.

\begin{lemma}\label{quasi.fr}
Let $A$ be an associative algebra over $k$, and let $A_{(k)}$,
$A_{(1)}=A$, $A_{(k)}=[A,A_{(k-1)}]$ be its central series with
respect to the commutator. If $A_{(p)}=0$ and $A_{(2)}^p=0$, then
setting $x \mapsto x^p$ defines an additive multiplicative central
map $\wt{\Fr}_k:A^\tw \to A$.
\end{lemma}

\proof{} Consider the tensor algebra $T^\hdot(A)$ generated by the
vector space $A$. Since $T^\hdot(A)$ is the universal enveloping
algebra for the free Lie algebra $L^\hdot(A)$, by the
Poincare-Birkhoff-Witt Theorem it has an increasing filtration
$F^{PBW}_\idot T^\hdot(A)$ such that the associated graded quotient
is isomorphic to the symmetric algebra $S^\hdot(L^\hdot(A))$.
Moreover, the universal enveloping algebra $T^\hdot(A)$ is a Hopf
algebra, this structure is compatible with the
Poincare-Birkhoff-Witt filtration, and it induces the standard Hopf
algebra structure on $S^\hdot(L^\hdot(A))$.

For any $x,y \in A$ let $d(x,y)=(x+y)^{\otimes p} - x^{\otimes p} -
y^{\otimes p} \in T^p(A)$. The element $d(x,y)$ is primitive with
respect to the Hopf algebra structure on $T^\hdot(A)$. {\em A
priori}, it lies inside the $p$-th term $F^{PBW}_p T^\hdot(A)$ of
the Poincare-Birkhoff-Witt filtration. However, the image of
$d(x,y)$ in $\gr^{PBW}_p T^p(A) \cong S^p(L^1(A)) = S^p(A)$ is
obviously zero, so that in fact $d(x,y) \in F^{PBW}_qT^\hdot(A)$ for
some $q \leq p-1$. Since $d(x,y)$ is primitive, so is its image in
$\gr^{PBW}_qT^\hdot(A) = S^q(L^\hdot(A))$; since $q < p$, this is
possible only if $q=1$. We conclude that $d(x,y) \in S^1(L^p(A)) =
L^p(A)$ is a Lie polynomial in $x$ and $y$.

But since by assumption $A_{(p)}=0$, the multiplication map
$m:T^\hdot(A) \to A$ vanishes on $L^p(A) \subset T^p(A)$. Therefore
$m(d(x,y)) = (x+y)^p-x^p-y^p = 0$, and $\wt{\Fr}_k$ is indeed
additive. Moreover, denoting $\ad x(z) = [x,z]$, we have
\begin{equation}\label{ad.fr}
[x^{\otimes p},y] = (\ad x)^p(y) \in L^{p+1}(A),
\end{equation}
so that $x^p$ is indeed central in $A$. Finally,
$$
\begin{aligned}
x^py^p &= \frac{1}{2}\left((x^p+y^p)^2-x^{2p}-y^{2p}\right) \\
&= \frac{1}{2}\left((x^2+y^2+2xy-[x,y])^p-x^{2p}-y^{2p}\right) =(xy)^p,
\end{aligned}
$$
which finishes the proof.
\endproof

Applying this to quantizations, we see that for any quantization
$A_h$ of a commutative algebra $A$, the quotient $A_h/h^{p-1}$
satisfies the conditions of the Lemma.

\begin{defn}\label{fr.const.defn} 
A {\em Frobenius-constant quantization} of a $k$-algebra $A$ is a
pair of a quantization $A_h$ of the algebra $A$ and a multiplicative
central $k[h]$-linear {\em splitting map} $s:A_h^\tw \to A_h$ such
that $s(h)=0$ and
$$
s(a) = a^p \mod h^{p-1}
$$
for any $a \in A_h$.
\end{defn}

Assume given a Frobenius-constant quantization $A_h$ of a
$k$-algebra $A$.  Since $s(h)=0$, the splitting map $s:A_h^\tw \to
A_h$ factors through the reduction map $A^\tw_h \to A^\tw \cong
A_h^\tw/h$. We will say that the Frobenius-constant quantization
$A_h$ is {\em regular} if $s$ further factors through the Frobenius
map $A^\tw \to A^p$ (in other words, $s(a)=0$ for any $a \in A^\tw$
with $a^p=0$). If the algebra $A$ is reduced, this is automatic; in
general, it might happen that for some $a \in A$, $a^p=0$ in $A$,
but when we lift $a$ to an element $\wt{a} \in A_h$, the $p$-th
power $\wt{a}^p$ becomes non-trivial already in $A_h/h^2$. We note
that in our approach to quantizations, we really have to consider
non-reduced algebras, even if the original algebra or scheme we want
to quantize is smooth (see Section~\ref{loc.sec}).

We also note that the notion of a regular Frobenius-constant
quantization has nice functoriality properties with respect to base
change. More precisely, assume given a $k$-algebra $A$, an ideal $I
\subset A^p$, and a regular Frobenius-constant quantization $\langle
A_h,s \rangle$ of $A$. Then the quotient $A_h/s(I)A_h$ is a regular
Frobenius-constant quantization of the quotient $A/IA$.

\subsection{Universal polynomials.}

To proceed further, we need to develop some formalism on the $p$-th
power operation in characteristic $p$. It will be convenient to
interpolate between the notions of an associative and a Poisson
algebra by introducing the following definition.

\begin{defn}
A {\em quantized algebra} $A$ over a field $k$ is an associative
algebra over the algebra $k[[h]]$ of formal power series in one
variable $h$ equipped with an additional $k[h]$-linear Lie bracket
$\{-,-\}$ which is a derivation in each variable and satisfies
\begin{equation}\label{bra.qu}
h\{a,b\}=ab-ba
\end{equation}
for every $a,b \in A$.
\end{defn}

A quantized algebra which is flat over $k[[h]]$ -- equivalently, has
no $h$-torsion -- is the same as a flat associative $k[[h]]$-algebra
$A$ with commutative quotient $A/hA$ (the bracket $\{a,b\}$ is
uniquely defined by \eqref{bra.qu}). On the other hand, a quantized
algebra annihilated by $h$ is the same as a Poisson algebra. If $V$
is finite-dimensional vector space over $k$, then the free
associative algebra over $k$ generated by $V$ is the tensor algebra
$T^\hdot V$. Since $T^\hdot(V)$ is also the universal enveloping
algebra of the free Lie algebra $L^\hdot(V)$ generated by $V$, it
has a Poincare-Birkhoff-Witt filtration $F_\idot^{PBW}T^\hdot V$,
and the free Poisson algebra $P^\hdot V$ is the associated graded
$\gr^{PBW}T^\hdot V$ with respect to the PBW filtration; as for the
free quantized algebra $Q^\hdot V$, it is the so-called {\em Reese
algebra} associated to the PBW filtration: we have
$$
Q^\hdot V = \bigoplus F_\idot T^\hdot V,
$$ 
and the variable $h$ acts as the natural embedding $h:F_\idot
T^\hdot V \to F_{\idot+1}T^\hdot V$. When $V=k \langle
x_0,\dots,x_n\rangle$ is spanned by basis elements $x_0,\dots,x_n$,
we will denote the algebra $Q^\hdot V$ by $Q^\hdot(x_0,\dots,x_n)$
and call its elements {\em quantized polynomials} in variables
$x_0,\dots,x_n$.

Consider now the algebra $Q^\hdot(x,y)$ of quantized polynomials in
two variables $x,y$. On one hand, this algebra has no $h$-torsion;
on the other hand, the quotient $Q^\hdot(x,y)/h^{p-1}$ obviously
satisfies the assumptions of Lemma~\ref{quasi.fr}. We conclude that
there exist two canonical quantized polynomials $L(x,y)$, $P(x,y)$
such that
\begin{equation}\label{poly}
\begin{aligned}
h^{p-1}L(x,y) &= (x+y)^p-x^p-y^p\\
h^{p-1}P(x,y) &= (xy)^p-x^py^p.
\end{aligned}
\end{equation}
As noted in the proof of Lemma~\ref{quasi.fr}, the polynomial
$L(x,y)$ is in fact a Lie polynomial in $x,y$. It is well-known, and
it can be computed explicitly (the Jacobson formula -- see
e.g. \cite[II, \S 7.3, D\'efinition 3.1]{dem2}). We will only need
an easy ``leading term'' formula; for the convenience of the reader,
we include a proof.

\begin{lemma}\label{sq}
Assume given a Lie algebra $L$ equipped with a two-step decreasing
filtration, $F^1L \subset F^0L=L$, $[F^i L,F^j L] \subset F^{i+j}L$
for any $i$, $j$. Then for any $x \in L$, $y \in F^1L$ we have
$$
L(x,y)=(\ad x)^{p-1}(y).
$$
\end{lemma}

\proof{} Since $\ad x$ preserves $F^1L$, and $F^2L=0$, so that the
bracket on $F^1L$ vanishes, the only non-trivial Lie monomial in $x$
and $y$ of arbitrary degree $l$ is $(\ad x)^l(y)$ (the proof is an
obvious induction on $l$). Therefore $L(x,y)=\lambda(\ad
x)^{p-1}(y)$ for some universal constant $\lambda \in
\Z/p\Z$. Consider now the Lie algebra $L$ spanned by $x$ and $y$,
with $F^1L \subset L$ spanned by $y$, and $\{x,y\}=y$. Then we have
\begin{align*}
\lambda y &= \lambda (\ad x)^p(y) = (\ad x)(\lambda (\ad x)^{p-1}(y))
= (\ad x)(L(x,y))\\ 
&= (\ad (x+y))^p(x) - (\ad x)^p(x) - (\ad y)^p(x) = y,
\end{align*}
which proves that $\lambda=1$.
\endproof

The Lie polynomial $L(x,y)$ is used to define the standard notion of
a so-called {\em restricted Lie algebra structure} in characteristic
$p$.

\begin{defn}
A {\em restricted Lie algebra} $A$ over a field of characteristic
$p$ is a Lie algebra $A$ over $k$ equipped with an additional
operation $x \mapsto x^{[p]}$ such that $(ax)^{[p]} = a^px^{[p]}$
for any $a \in k$, $x \in A$, and
\begin{align}
\{x^{[p]},y\} &= (\ad x)^p(y),\label{restr.lie.1}\\
(x+y)^{[p]} &= x^{[p]} + y^{[p]} + L(x,y)\label{restr.lie.2}
\end{align}
for any $x,y \in A$, where we denote $\ad x(z)=\{x,z\}$.
\end{defn}

We will need a Poisson version of this notion.

\begin{defn}\label{restr.poi.def}
A {\em restricted Poisson algebra} $A$ over a field $k$ of
characteristic $p$ is a Poisson algebra $A$ over $k$ equipped with an
operation $x \mapsto x^{[p]}$ which turns $\langle A,
\{-,-\}\rangle$ into a restricted Lie algebra, and satisfies
\begin{equation}\label{restr.poi}
(xy)^{[p]}  = x^py^{[p]} + x^{[p]}y^p + P(x,y)
\end{equation}
for any $x,y \in A$.
\end{defn}

More generally, here is a version for arbitrary quantized algebras.

\begin{defn}\label{restr.qua.def}
A {\em restricted structure} on a quantized algebra $A$ over $k$ is
given by an operation $x \mapsto x^{[p]}$ on $A$ which turns
$\langle A, \{-,-\}\rangle$ into a restricted Lie algebra, preserves
$h$ -- that is, $h^{[p]}=h$ -- and satisfies
$$ 
(xy)^{[p]} = x^py^{[p]} + x^{[p]}y^p - h^{p-1}x^{[p]}y^{[p]} +
P(x,y)
$$ 
for any $x,y \in A$.
\end{defn}

As an immediate corollary of this definition, we see that for any
element $x \in A$ of a restricted quantized algebra $A$ we have
\begin{equation}\label{restr.h}
(hx)^{[p]}=hx^p.  
\end{equation}
If $A$ is a Poisson algebra, we have $hA=0$, so that these two
definitions agree. If the Poisson bracket on $A$ is trivial, then
the $p$-th power operation is given by a Frobenius-derivation $K:A
\to A$; thus a restricted Poisson algebra with trivial bracket is
the same as a commutative algebra $A$ equipped with a
Frobenius-derivation $K:A \to A$. Even if the bracket on $A$ is not
trivial, the {\em difference}\label{dff} between two restricted
structures is a Frobenius-derivation of $A$ (moreover, by
\eqref{restr.lie.1} it must take values in the Poisson center of
$A$). In general, examples of restricted quantized algebras are
provided by Frobenius-constant quantizations. Namely, assume given a
Frobenius-constant quantization $A_h$ of a commutative algebra $A$
over $k$, with splitting map $s:A_h^\tw \to A_h$. Then for any
integer $n$, the quotient $A_n=\wt{A}/h^{n+1}$ is a quantized
algebra, and this algebra is restricted: a natural $p$-th power
operation is given by
$$ 
x^{[p]} = \frac{1}{h^{p-1}}(x^p - s(x)).
$$
Checking the conditions of Definition~\ref{restr.qua.def} is
straightforward and left to the reader. Conversely, every restricted
quantized algebra $A_h$ without $h$-torsion is a Frobenius-constant
quantization of its quotient $A=A_h/h$: using \eqref{restr.lie.1},
\eqref{restr.lie.2}, and \eqref{restr.h}, one easily checks that
$s(x)=x^p - h^{p-1}x^{[p]}$ gives a central splitting map $s:A_h^\tw
\to A_h$ which vanishes on $hA_h^\tw \subset A_h^\tw$.

\subsection{Restricted structures.}\label{prp.sub}
We can now state our results. We fix the base field $k$ of odd
positive characteristic, and we only consider algebras and schemes
over $k$. We will say that a Noetherian scheme $X$ equipped with a
map $\pi:X \to S$ to a Noetherian scheme $S$ is {\em quasiregular}
over $S$ if $\pi$ is flat, $X$ is of finite type over $S$, and the
sheaf $\Omega^1(X/S)$ of relative K\"ahler differentials is flat
over $\calo_X$. A scheme $X$ is {\em quasiregular} if it is
quasiregular over the point $\Spec k$. A regular scheme $X$ of
finite type is of course quasiregular; since $p=\cchar k$ is
positive, the converse is true only if $X$ is reduced (and at least
one non-reduced quasiregular algebra will be really important to us,
see Section~\ref{loc.sec} below). Most of the usual facts of
differential calculus works for quasiregular schemes without any
changes -- in particular, we have flat sheaves $\Omega^\hdot(X/S)$
of relative differential forms, a flat relative tangent bundle
$\T(X/S)$, the de Rham differential on $\Omega^\hdot(X/S)$, and the
Cartan homotopy formula describing the natural action of the Lie
algebra $\T(X/S)$ on $\Omega^\hdot(X/S)$.

It is customary to consider the pullback $X \times_{\Fr_S} S$ of
$X/S$ with respect to the Frobenius map $\Fr_S:S \to S$, and to
factor the Frobenius map $\Fr_X:X \to X$ through the pullback map
$\id \times \Fr_S:X \times_{\Fr_S} S \to X$ by the relative
Frobenius map $\Fr_{X/S}:X \to X \times_{\Fr_S} S$. In other words,
one considers the commutative diagram
$$
\begin{CD}
X @>{\Fr_{X/S}}>> X \times_{\Fr_S} S @>{\id \times \Fr_S}>>
X\\
@. @VVV @VV{\pi^\tw}V\\
@. S @>{\Fr_S}>> S
\end{CD},
$$
where the square on the right is Cartesian. However, it will be more
convenient for us to consider the sheaf $\calo_{X}^{[p]} = \calo_X^p
\otimes_{\calo_{S}^p} \calo_S$ of algebras on $X$ and the scheme
$X^{[p]} = \langle X, \calo_{X}^{[p]}\rangle$. If $X$ and $S$ are
reduced, we have $X^{[p]} \cong X \times_{\Fr_S} S$. If if they are
not reduced, the scheme $X^{[p]}$ is still quasiregular, and the
Frobenius map $\Fr_X:X \to X$ still factors through a map
$\Fr_{X^{[p]}}: X^{[p]} \to X$ by means of a map which we will, by
abuse of notation, denote by $\Fr_{X/S}:X \to X^{[p]}$.

One easily checks that the sheaf $\calo_X^{[p]}$ of functions on
$X^{[p]}$ coincides with the subsheaf $\calo_{X,cl} \subset \calo_X$
of functions closed with respect to the relative de Rham
differential. Therefore the de Rham differential on $X/S$ is
$\calo_X^{[p]}$-linear, so that the de Rham complex
$\Fr_{X/S*}\Omega^\hdot(X/S)$ is a complex of coherent sheaves of
$\calo_X^{[p]}$-modules on the scheme $X^{[p]}$. In particular, we
have the coherent $\calo_X^{[p]}$-module subsheaves
$\Fr_{X/S*}\Omega^\hdot_{cl}(X/S) \subset
\Fr_{X/S*}\Omega^\hdot(X/S)$ of closed forms.

The Zariski cohomology sheaves $\hh^\hdot(\Fr_{X/S*}\Omega^\hdot
(X/S))$ are identified with the sheaves $\Omega^\hdot_{[p]}(X/S) =
\Fr_{X^{[p]}}^*\Omega^\hdot(X/S)$ by the canonical
$\calo_X^{[p]}$-linear Cartier operations
$C^\hdot:\Fr_{X/S*}\Omega^\hdot_{cl}(X/S) \to
\Omega^\hdot_{[p]}(X/S)$ (this generalizes the isomorphism
$\calo_{X,cl} \cong \calo_X^{[p]}$). The canonical filtration on
$\Omega^\hdot(X/S)$ induces an increasing filtration $F_\idot$ on
the relative de Rham cohomology sheaves
$\RR^\hdot\pi_*(\Omega^\hdot(X/S)$ known as the {\em conjugate
filtration}; the Cartier operations $C^\hdot$ give the projections
$F_j\RR^k\pi_*(\Omega^\hdot(X/S)) \to
\RR^{k-j}\pi_*\Omega^j_{[p]}(X/S)$ onto the associated graded pieces
of this filtration.

We also apply the standard deformation theory to the family $X/S$;
thus for any derivation $K:\calo_S \to \calo_S$, we have the
Kodaira-Spencer deformation class $\theta_K \in
H^1(X,\T(X/S))$. This class is represented by an explicit
$\T(X/S)$-torsor $[\theta_K]$ over $X$: local sections of
$[\theta_K]$ are derivations of the structure sheaf $\calo_X$ whose
restrictions to $\pi^{-1}\calo_S \subset \calo_X$ coincide with $K$.
For any relative $k$-form $\alpha \in H^0(X,\Omega^k(X/S))$,
contraction with $\alpha$ gives a sheaf map $\T(X/S) \to
\Omega^{k-1}(X/S)$, $\xi \mapsto \alpha \cntrct \xi$. By means of
this map, the torsor $[\theta_K]$ induces an
$\Omega^{k-1}(X/S)$-torsor on $X$ which we will denote by $\alpha
\cntrct [\theta_K]$. This immediately generalizes to derivations $K$
with values in some sheaf of $\calo_S$-modules; in particular, for
any Frobenius-derivation $K:\calo_S \to \calo_S$ and any form
$\alpha \in \Omega^k(X/S)$, we obtain an
$\Omega^{k-1}_{[p]}(X/S)$-torsor $\Fr_{X^{[p]}}^*\alpha \cntrct
[\theta_K]$.

Finally, for any vector field $\xi$ on $X$, the $p$-th power of the
corresponding derivation of the structure sheaf $\calo_X$ is also a
derivation; therefore the Lie algebra $\T(X/S)$ of relative vector
fields on $X$ carries a natural restricted Lie algebra structure.

By a symplectic manifold over $S$ we will understand a Noetherian
quasi\-regular scheme $X/S$ equipped with a non-degenerate closed
relative $2$-form $\Omega \in H^0(X,\Omega^2(X/S))$. Just as in the
case of characteristic $0$, such a manifold is automatically equipped
with an $\calo_S$-linear Poisson bracket. For any local function $f
\in \calo_X$ we denote the derivation $g \mapsto \{f,g\}$ of the
sheaf $\calo_X$ by $H_f$; explicitly, we have $H_f \cntrct \Omega =
df$. Vector fields $H_f$ are called Hamiltonian and form a subsheaf
of Lie subalgebras $H(X) \subset \T(X/S)$ in the relative tangent
sheaf $\T(X/S)$.

\begin{lemma}\label{no.poi.cent}
Assume given a symplectic manifold $X/S$. Then the Poisson center
$\ZZ \subset \calo_X$ coincides with the subsheaf $\calo_X^{[p]}
\subset \calo_X$.
\end{lemma}

\proof{} A function $f \in \calo_X$ lies in the Poisson center if
and only if $H_f=0$. But since $\Omega$ is non-degenerate, this is
equivalent to $H_f \cntrct \Omega =0$, which means $df=0$.
\endproof

As a corollary of this Lemma, we see that for symplectic manifolds
over $k$, a quantization $\calo_h$ is Frobenius-constant if and only
if \thetag{a} it is central, and \thetag{b} the center $\ZZ_h$ is
isomorphic to $\calo_X^p[[h]]$ -- indeed, any splitting map
$\calo_X^p \cong \ZZ \to \ZZ_h$ induces an isomorphism
$\calo_X^p[[h]] \cong \ZZ_h$ and conversely, any map $\calo_X^p \to
\calo_X^p[h]/h^{p-1} \cong \ZZ_h/h^{p-1}$ can be obviously lifted to
a map $\calo_X^p \to \calo_X^p[[h]]$. Therefore -- at least in the
case $S = \Spec k$ -- our Definition~\ref{fr.const.defn} and
\cite[Definition 3.3]{BK2} impose equivalent conditions on
quantizations. The difference is in the additional structure (not
every splitting map in \cite[Definition 3.3]{BK2} would serve for
Definition~\ref{fr.const.defn}).

Our first general result in this paper is the following.

\begin{theorem}\label{cent}
Assume given a symplectic manifold $X/S$. Then the following
conditions are equivalent:
\begin{enumerate}
\item\label{ham.restr} The subalgebra of Hamiltonian vector fields
in $\T(X/S)$ is closed with respect to the restricted Lie algebra
operation $\xi \mapsto \xi^{[p]}$.
\item We have $C^2(\Omega)=0 \in H^0(X,\Omega^2_{[p]}(X/S))$.
\end{enumerate}
Moreover, if $X$ admits a central quantization $\calo_h$, then both
these conditions are satisfied.
\end{theorem}

The proof of this Theorem is given in Subsection~\ref{cent.sub}
(page \pageref{cent.pf}). As an immediate corollary, we see if a
symplectic manifold $X$ admits a central quantization, or at least
satisfies the condition \ref{ham.restr} of Theorem~\ref{cent}, then
the de Rham cohomology class $[\Omega] \in H^2_{DR}(X/S)$ lies in
the first term $F_1H^2_{DR}(X/S)$ of the conjugate filtration. In
fact, $[\Omega]$ comes from a cohomology class
$$
\overline{[\Omega]} \in H^1(X,\Omega^1_{cl}(X/S)).
$$ 
Indeed, geometrically Theorem~\ref{cent}~\ref{ham.restr} means that
the form $\Omega$ is locally exact in Zariski topology. The
$1$-forms $\alpha$ such that $d\alpha = \Omega$ form a torsor
$\M_\Omega$ over the vector bundle $\Fr_{X/S*}\Omega^1_{cl}(X/S)$ on
$X^{[p]}$ spanned by closed $1$-forms on $X$; this torsor represents
the class $\overline{[\Omega]}$. Applying the Cartier map
$C^1:\Fr_{X/S*}\Omega^1_{cl}(X/S) \to \Omega^1_{[p]}(X/S)$, we can
also consider the induced $\Omega^1_{[p]}(X/S)$-torsor
$C^1_*\M_\Omega$.

Assume now that the base scheme $S$ is equipped with a
Frobenius-de\-ri\-va\-ti\-on $K:\calo_S \to \calo_S$, so that
$\calo_S$ becomes a restricted Poisson algebra with trivial Poisson
bracket. By a restricted structure on the Poisson scheme $X/S$ we
will understand a restricted structure on the Poisson algebra sheaf
$\calo_X$ such that the embedding $\pi^{-1}\calo_S \to \calo_X$ is a
(central) restricted Poisson algebra map. Then we have the following
(the proof is in Subsection~\ref{obstr.sub}, on page
\pageref{obstr.pf}).

\begin{theorem}\label{obstr}
Let $X$ be a symplectic manifold over $S$ which satisfies the
equivalent conditions of Theorem~\ref{cent}. Then restricted
structures on the Poisson scheme $X/S$ are in natural one-to-one
correspondence with isomorphisms
$$
C^1_*\M_\Omega \cong \Fr_{X^{[p]}}^*\Omega \cntrct [\theta_K]
$$
of torsors over the bundle $\Omega^1_{[p]}(X/S)$.
\end{theorem}

This Theorem immediately gives a necessary condition for the
existence of Frobenius-constant quantizations. Namely, assume given
a symplectic manifold $X/S$, $\pi:X \to S$, and a
Frobenius-derivation $K:\calo_S \to \calo_S$. We will say that a
regular Frobenius-constant quantization $\langle \calo_h,s \rangle$
of the Poisson scheme $X$ is compatible with the derivation $K$ if
for any local section $f \in H^0(U,\calo_h)$ whose reduction
$\overline{f} \in H^0(U,\calo_X)$ modulo $h$ comes from a function
$b \in H^0(S,\calo_S)$ on $S$, $\overline{f}=\pi^*(b)$, we have
\begin{equation}\label{good.on.base}
f^p-s(f^p)=\pi^*(K(b)).
\end{equation}

\begin{corr}
Assume given a symplectic manifold $X/S$, a Frobenius-derivation
$K:\calo_S \to \calo_S$, and a regular Frobenius-constant
quantization $\langle \calo_h,s \rangle$ of the Poisson scheme $X$
compatible with $K$. Then the cohomology class $[\Omega] \in
H^2_{DR}(X/S)$ of the symplectic form $\Omega$ lies in the first
term $F_1H^2_{DR}(X/S)$ of the conjugate filtration, and satisfies
\begin{equation}\label{G.M.C}
C^1([\Omega]) = (\Fr_{X^{[p]}}^*[\Omega]) \cntrct \theta_K 
\in H^1(X,\Omega^1_{[p]}(X/S)).
\end{equation}
In particular, if $K=0$ -- for instance, if $S=\Spec k$ -- then the
class $[\Omega]$ lies in the second term $F_2H^2_{DR}(X/S)$ of the
conjugate filtration.
\end{corr}

\proof{} As we have noted in the definition of the restricted
quantized algebras, setting 
$$
x^{[p]} = x^p - \frac{1}{h^{p-1}}s(x)
$$
gives a restricted quantized algebra structure on the sheaf
$\calo_h$. Reducing modulo $h$, we obtain a restricted Poisson
structure on the quotient $\calo_h/h \cong \calo_X$. Compatibility
condition \eqref{good.on.base} insures that this structure is
compatible with the restricted Poisson structure on $\calo_S$ given
by $K$.  We can now apply Theorem~\ref{cent}~\ref{ham.restr} and
Theorem~\ref{obstr}, and the only thing that remains to be proved is
that the quantization $\calo_h$ is central. Indeed, the splitting
map $s$ induces a map $s:\calo_X^p \otimes_{B^p} B_h \to \ZZ_h$,
which reduces to $\calo_X^p \otimes_{B^p} B \to \ZZ_h/h \subset \ZZ$
modulo $h$. But by Lemma~\ref{no.poi.cent} the composition
$\calo_X^p \otimes_{B^p} B \cong \calo_{X/S}^p \to \ZZ$ is an
isomorphism. Therefore $\ZZ_h/h \cong \ZZ$.
\endproof

\begin{remark}\label{proj.rem}
The conditions $C([\Omega]) = 0$ and furthermore $C^1([\Omega])=0$
are not very natural for {\em projective} symplectic manifolds such
as surfaces of type $K3$ and abelian varieties of even dimension. At
least for abelian varieties, a more natural condition is
$C([\Omega])=[\Omega]$. It would be very interesting to develop a
parallel theory for symplectic manifolds of this type. We do not
know what should replace the condition $C^1([\Omega])=0$.
\end{remark}

\subsection{Quantizations.}\label{exi.sub}
We now turn to existence results for quantizations. We change our
point of view: we assume from now on that the base scheme $S$ is
local, we replace $X/S$ with its special fiber, and we incorporate
both quantizations and the family $X/S$ into a single
multi-parameter partially non-commutative deformation of the
structure sheaf. To package the data, we introduce the following.

\begin{defn}\label{good.base.defn}
A {\em quantization base} is a commutative Noetherian complete local
$k[[h]]$-algebra $B$ with maximal ideal $\m_B \subset B$, $h \in
\m_B$, $k \cong B/\m_B$ equipped with a map $K:B \to B^\tw$ such
that setting $b^{[p]} = K(b)$ for any $b \in B$ turns $B$ into a
restricted quantized algebra (with trivial bracket).
\end{defn}

Assume given a symplectic manifold $X/k$, and a quantization base
$B$ with maximal ideal $\m_B \subset B$. By a {\em $B$-quantization}
$\calo^B$ of the manifold $X$ we will understand a sheaf $\calo^B$
of restricted quantized flat $B$-algebras on $X$ equipped with a
restricted Poisson isomorphism isomorphism $\calo^B/\m_B\calo^B
\cong \calo_X$. In other words, $\calo^B$ is equipped with a central
restricted quantized map $B \to H^0(X,\calo^B)$, $\calo^B$ with the
induced $B$-module structure is flat over $B$, and the quotient
$\calo^B/\m_B\calo^B$ is identified with $\calo_X$. In the case $B =
k[[h]]$, this reduces to Definition~\ref{fr.const.defn}, so that a
$k[[h]]$-quantization is the same as a Frobenius-constant
quantization. Another basic example is $B = k$; a $k$-quantization
is the same as a restricted structure on the Poisson scheme $X$. We
will say that a $B$-quantization $\calo^B$ is {\em regular} if for
any local section $f \in H^0(U,\calo^B)$ whose reduction
$\overline{f} \in H^0(U,\calo_X)$ mod $\m_B$ satisfies
$\overline{f}^p=0$, we have $f^p=h^{p-1}f^{[p]}$. As before, this
condition is automatic if $X/k$ is reduced.

\medskip

We construct quantizations by induction, going step-by-step through
a sequence of extensions of the quantization base. To make this
precise, we say that by an {\em extension} $\langle B, I \subset
B\rangle$ of a quantization base $B_0$ we will understand a
quantization base $B$ equipped with a (restricted) ideal $I \subset
B$ and an isomorphism $B/I \cong B_0$. We will say that an extension
$\langle B,I \rangle$ is {\em small} if the ideal $I \subset B$ is
annihilated by $\m_B$. For example, for any quantization base $B$
and any integer $n$, the quotient $B/\m_B^{n+1}$ is a small
extension of the quotient $B/\m_B^n$.

Our main result gives the induction step -- we describe in detail
how to prolong a quantization over a small extension of its
base. Namely, by a {\em small Dieudonn\'e module} over $k$ we will
understand a vector space $I$ equipped with a $p$-linear map $K_I:I
\to I$, $K_I(av)=a^pK_I(v)$ for any $a \in k$, $v \in I$ (for a
reader familiar with the notion of a Dieudonn\'e module, we note
that ``small'' here means that the map $V$ is trivial, $V=0$, and we
take $F=K$). We note that for any small extension $\langle B,I
\rangle$, the ideal $I$ carries a natural structure of a small
Dieudonn\'e module: we set $K_I(x) = x^{[p]}$.

\begin{defn}
Let $X/k$ be a symplectic manifold, and assume given a small
Dieudonn\'e module $I$ over $k$. By
$$
\hh\langle I \rangle \subset \Fr_*\Omega^1_{cl} \otimes_k I
$$
we will denote the sheaf of closed $I$-valued $1$-forms $\alpha$ on
$X$ satisfying $C^1(\alpha) = K_I(\alpha)$.
\end{defn}

\begin{remark}\label{not.lin}
We note that unless $K_I=0$, the sheaf $\hh\langle I \rangle$ is
{\em not} a sheaf of $\calo_X^{[p]}$-modules. Indeed, the Cartier
map $C^1$ is $\calo_X^{[p]}$-linear, but the map $K_I$ is not: we
have $K_I(f\alpha)=f^pK_I(\alpha)$.
\end{remark}

\begin{prop}\label{step}
Assume given a small extension $\langle B,I \rangle$ of a
quantization base $B_0$, and a regular $B_0$-quantization $\calo_h$
of a symplectic manifold $X/k$. Denote by $Q(B,\calo_h)$ the set of
isomorphism classes of all $B$-quantizations $\calo_h'$ of $X/k$
equipped with an isomorphism $\calo_h'/I \cong \calo_h$. Then the set
$Q(B,\calo_h)$ is described by the \'etale cohomology of the sheaf
$\hh\langle I \rangle$ on $X$ in the following way:
\begin{enumerate}
\item There exists a canonical obstruction class $c(\calo_h) \in
  H^2_{et}(X,\hh \langle I \rangle)$ such that $Q(B,\calo_h)$ is
  non-empty if and only if $c(\calo_h)=0$.
\item If the set $Q(B,\calo_h)$ is not empty, then it has a natural
  structure of a torsor over the group $H^1_{et}(X,\hh\langle I
  \rangle)$.
\end{enumerate}
\end{prop}

This Proposition is proved in Subsection~\ref{step.pf.sub} on
page~\pageref{step.pf}. As usual, one can give a more precise result
which also describes the category of the quantizations $\calo_h'$
and its gerb structure; we do not do this to save space. 

Our second result is somewhat surprising.

\begin{prop}\label{level.1}
In the assumptions of Proposition~\ref{step}, assume in addition
that $B_0=k$, so that $I = \m_B$ and $B$ is an arbitrary
quantization base with $\m_B^2=0$. Assume also given a restricted
structure on the Poisson scheme $X/k$, so that $\calo_X$ becomes a
$k$-quantization.  Then the obstruction class $c(\calo_X)$
tautologically vanishes, and the torsor $Q(B,\calo_X)$ is trivial:
we have a natural identification
$$
Q(B,\calo_X) \cong H^1(X,\hh\langle \m_B \rangle).
$$
\end{prop}

The proof is also in Subsection~\ref{step.pf.sub}, on
page~\pageref{step.pf}. To get a feeling for what
Proposition~\ref{step} and Proposition~\ref{level.1} mean in
practice, we need to get a handle on the sheaves $\hh \langle V
\rangle$. Specifically, we need one particular case: $V=k$ is
one-dimensional, and $K:V \to V$ is given by $K(a)=a^p\lambda$ for
some constant $\lambda \in k$. We denote this Dieudonn\'e module by
$k\langle \lambda \rangle$, and we denote the corresponding sheaf by
$\hh\langle\lambda\rangle$. If $\lambda=0$, then $\hh = \hh\langle 0
\rangle$ is the sheaf of closed $1$-forms $\alpha$ on $X$ with
$C^1(\alpha)=0$ -- in other words, the sheaf of exact
$1$-forms. Thus we have a short exact sequence
\begin{equation}\label{seq.1}
\begin{CD}
0 @>>> \calo_X^p @>>> \calo_X @>>> \hh @>>> 0.
\end{CD}
\end{equation}
The right-hand side map is the de Rham differential, and the
left-hand side map is the Frobenius map on $X$. If $\lambda = 1$,
then $\hh_{log}=\hh \langle 1 \rangle$ is the sheaf of closed
$1$-forms with $C^1(\alpha)=\alpha$. It is well-known that these are
of the form $\alpha = df/f$, where $f$ is an invertible function on
$X$. Thus we have a short exact sequence
\begin{equation}\label{seq.2}
\begin{CD}
0 @>>> \left(\calo_X^p\right)^* @>>> \calo_X^* @>>> \hh_{log} @>>>
0.
\end{CD}
\end{equation}
This is a version of the Kummer exact sequence. The map on the
right-hand side is the $d\log$-map, $f \mapsto df/f$, and the map on
the left-hand side is again the Frobenius, which for the
multiplicative group $\calo_X^*$ means just the multiplication by
$p$.

\begin{remark}
Both \eqref{seq.1} and \eqref{seq.2} are examples of a more general
construction. The category of small Dieudonn\'e modules is
equivalent to the category of algebraic groups over $k$ annihilated
by $p$ (see e.g. \cite{Dem}). Algebraic groups over $k$ are the same
as groups sheaves on $\Spec k$ in flat topology; by pullback, an
algebraic group $G$ defines a sheaf $\GG$ on $X$ in flat
topology. If we denote by $\rho$ the natural map from $X$ with flat
topology to $X$ with \'etale topology, then $\RR^i\rho_*\GG=0$
unless $i=1$. And if $G=G_I$ corresponds to a Dieudonn\'e module
$V$, then $\RR^1\rho_*\GG \cong \hh\langle I \rangle$. The sequence
\eqref{seq.1} computes this higher direct image for $k\langle 0
\rangle$ (which corresponds to the group of points of order $p$ in
$\ga$), while \eqref{seq.2} corresponds to $k\langle 1 \rangle$ (the
group of points of order $p$ in $\gm$). The reader can find more
details, for instance, in \cite[III, \S 4]{mln}.
\end{remark}

Unfortunately, we cannot say anything about the obstruction class
$c(\calo_h)$ in Proposition~\ref{step}~\thetag{i}; to get results,
we have to kill it off by imposing a condition on the symplectic
manifold $X$.

\begin{defn}
A scheme $X$ is called {\em admissible} if it is reduced, and the
Frobenius map $\Fr:H^i(X^{[p]},\calo_X^{[p]}) \to H^i(X,\calo_X)$
is bijective for $i=1,2,3$.
\end{defn}

By the cohomology long exact sequence associated to \eqref{seq.1},
we see that an admissible scheme $X$ satisfies
$$
H^1(X,\hh) \cong H^2(X,\hh)=0.
$$
The sheaf $\hh$ is in fact a coherent sheaf of $\calo_X^p$-modules,
so its cohomology is the same in all topologies one might want to
use, including the Zariski topology.\ We note that the admissibility
condition always holds if $X$ is reduced, projective and ordinary
over $S = \Spec k$, or if for some reason $H^i(X,\calo_X) = 0$ for
$i=1,2,3$ (by the Grauert-Riemenschneider Theorem, this is always
the case when $X$ is projective over an affine $Y$ and lifts to a
smooth manifold over the ring $W_2(k)$ of second Witt vectors of the
field $k$).

We can now deduce existence results for quantizations of admissible
symplectic manifolds.

\begin{prop}\label{B.qua}
Let $X/k$ be an admissible symplectic manifold equipped with a
restricted structure, and let $B$ be a quantization base with
maximal ideal $\m_B \subset B$. Then the set $Q(B,X)$ of isomorphism
classes of regular $B$-quantizations $\calo^B$ of $X$ equipped with a
restricted Poisson isomorphism $\calo^B/\m_B \cong \calo_X$ is
naturally identified with the cohomology group
$H^1_{et}(X,\hh\langle\m_B/\m_B^2\rangle)$.
\end{prop}

\proof{} Filter the local algebra $B$ by powers of the maximal ideal
$\m_B$, and denote $B_n = B/\m_B^{n+1}$. Then $B_{n+1}$ is a small
extension of $B_n$, with ideal $I = \m_B^{n+1}/\m_B^{n+2}$. However,
it immediately follows from Definition~\ref{restr.qua.def} that
whenever $n \geq 1$, for every $b \in \m_B^{n+1} \subset B$ we have
$b^{[p]} \in \m_B^{n+2} \subset B$. Therefore the map $K:I \to I$
given by the restricted structure is trivial for $n \geq 1$. By
Proposition~\ref{step}, we conclude that every regular
$B_1$-quantization of $X$ extends uniquely to a $B$-quantization. It
remains to notice that $B_1 = B/\m_B^2$ satisfies the assumptions of
Proposition~\ref{level.1}.
\endproof

\begin{theorem}\label{qua}
Assume given a quantization base $B$ and a symplectic ma\-n\-ifold
$X/S$ over $S = \Spec (B/hB)$ such that the special fiber $X_o =
X/\m_B$ is admissible. Assume also given a restricted structure on
the Poisson scheme $X/S$ compatible with the natural restricted
structure on $B/hB$. Then this restricted structure extends to a
regular Frobenius-constant quantization of $X$ compatible with the
restricted structure on $B$ as in \eqref{good.on.base}. Moreover,
the set $Q(B,\calo_X)$ of isomorphism classes of all such
quantizations is naturally a torsor over the group
$H^1_{et}(X,\hh_{log})$.
\end{theorem}

\proof{} Consider the sets $Q(B,X_o)$, $Q((B/hB),X_o)$ of
isomorphism classes of all regular $B$-quantizations,
resp. $(B/hB)$-quantization of the restricted Poisson scheme
$X_o$. Equivalently, $Q((B/hB),X_o)$ is the set of isomorphism
classes of restricted Poisson deformations $X/S$ of $X_o$. Reduction
modulo $h$ defines a map $\sigma:Q(B,X_o) \to Q((B/hB),X_o)$.

The given family $X/S$ defines a point $p \in Q((B/hB),X_o)$, and we
have to prove that the preimage $\sigma^{-1}(p) \subset Q(B,X_o)$ is
not empty and is a torsor over $H^1(X,\hh_{log})$. Indeed, denote $I
= \m_B/\m_B^2$ and $I_0 = \m_{B/hB}/\m_{B/hB}^2$. We have a
short exact sequence
\begin{equation}\label{sp.seq}
\begin{CD}
0 @>>> h \cdot k @>>> I @>>> I_0 @>>> 0
\end{CD}
\end{equation}
of small Dieudonn\'e modules, and Definition~\ref{restr.qua.def}
implies in particular that $K(h) = h^{[p]}=h$, so that the module on
the left-hand side is $k\langle 1 \rangle$. Therefore we have a
short exact sequence of sheaves
\begin{equation}\label{sh.seq}
\begin{CD}
0 @>>> \hh_{log} @>>> \hh\langle I \rangle @>>> \hh\langle I_0
\rangle @>>> 0
\end{CD}
\end{equation}
on $X$, and by Proposition~\ref{B.qua}, $Q(B,X_o)$ and
$X((B/hB),X_o)$ are naturally isomorphic to the first \'etale
cohomology group of $X$ with coefficients in the second and the
third term of \eqref{sh.seq}. Thus to prove the claim, it suffices
to show that taking first cohomology preserves the exactness of
\eqref{sh.seq}. For this it suffices to prove that \eqref{sh.seq}
splits, and it further suffices to prove that \eqref{sp.seq}
splits. This is well-known: $k\langle 1 \rangle$ is an injective
object in the category of small Dieudonn\'e modules over $k$ (see
e.g. \cite{dem2}).
\endproof

Theorem~\ref{qua} is analogous to \cite[Theorem 1.8]{BK} in the
characteristic $0$ case. We note that in characteristic $0$, one has
to consider information on all levels to obtain a full
parameterization of quantizations. In positive characteristic, we
cannot really do it; we restrict our attention to Frobenius-constant
quantizations, impose a stronger admissibility condition, and
extension to higher levels becomes automatic and
unique. Unfortunately, if a symplectic manifold $X$ satisfying
assumptions of Theorem~\ref{qua} is projective over $k$, then we
must have $[\Omega]=0$, which is impossible. Therefore, as we have
already noted in Remark~\ref{proj.rem}, our methods do not apply
well enough to projective symplectic manifolds such as surfaces of
type $K3$. Most probably, such manifolds do not admit
Frobenius-constant quantizations at all.

Our final result concerns the particular case of Theorem~\ref{qua}
when $B=k[[h]]$, $S = \Spec k$, so that $X=X_o$ is itself admissible
(in particular, it is reduced, and we have $X^{[p]}=X^\tw$). In this
case, the torsor in Theorem~\ref{qua} is obviously split, so that we
have a natural isomorphism
$$
Q(k[[h]],X) \cong H^1_{et}(X,\hh_{log}).
$$
In other words, every regular Frobenius-constant quantization
$\calo_h$ of $X$ is uniquely characterized by a cohomology class
$[\calo_h] \in H^1(X,\hh_{log})$. The cohomology long exact sequence
associated to \eqref{seq.2} allows to give a more precise
description of the group $H^1_{et}(X,\hh_{log})$: we have a short exact
sequence
$$
\begin{CD}
0 @>>> \Pic(X)/p\Pic(X) @>>> H^1_{et}(X,\hh_{log}) @>{\Br}>>
\Br_p(X) @>>> 0,
\end{CD}
$$
where $\Pic(X) = H^1_{et}(X,\calo_X^*)$ is the Picard group of $X$,
and $\Br_p(X) \subset H^2_{et}(X,\calo_X^*)$ is the $p$-torsion part
of the Brauer group. By projection, the characterizing class
$[\calo_h] \in H^1_{et}(X,\hh_{log})$ of every regular
Frobenius-constant quantization $\calo_h$ of $X$ gives a $p$-torsion
element $\Br([\calo_h])$ in the Brauer group.

This has the following interpretation in terms of
quantizations. Consider the scheme $\wh{X} = \sspec \calo_X^p[[h]]$;
this is essentially the completion of $X^\tw \times \Spec k[[h]]$
along the special fiber $X^\tw \subset X^\tw \times \Spec k[[h]]$,
but it exists as an actual scheme over $\Spec k[[h]]$ rather than
simply a formal scheme, so that the complement $\overline{X} =
\wh{X} \setminus X^\tw$ is a well-defined scheme as well. Denote by
$\rho:\overline{X} \to X^\tw$ the projection. By \cite[III,
Th\'eor\`eme 5.4.5]{EGA}, the category of coherent sheaves on
$\wh{X}$ is equivalent to the category of complete finitely
generated topological sheaves of $\calo_X^p[[h]]$-modules on
$X^\tw$. In particular, any Frobenius-constant quantization
$\calo_h$ of $X$ defines a coherent sheaf of algebras on $\wh{X}$.

\begin{prop}\label{azu}
For any central quantization $\calo_h$ of a symplectic manifold
$X/S$, the algebra $\calo_h(h^{-1})$ is a simple algebra over the
localization $\ZZ_h(h^{-1})$ of its center $\ZZ_h$. Moreover, assume
that $S=\Spec k$ and that we have $H^2(X,\calo_X)=0$. Then for any
regular Frobenius constant quantization $\calo_h$ of $X$, the class
of the Azumaya algebra $\calo_h(h^{-1})$ in the Brauer group
$\Br(\overline{X})$ is equal to the pullback $\rho^*\Br([\calo_h])$
of the projection $\Br([\calo_h])$ of the characterising class
$[\calo_h] \in H^1_{et}(X,\hh_{log})$.
\end{prop}

This is proved in Subsection~\ref{azu.sub} on
page~\pageref{azu.pf}. In particular, this means that an admissible
$X/k$ has quantizations which give a split Azumaya algebra on
$\overline{X}$. 

\subsection{Conjectures.}
We would like to finish the section with some speculations motivated
by Lemma~\ref{quasi.fr}. Namely, we propose the following.

\begin{conj}\label{cnj}
Let $B$ be a local commutative algebra over $k$ with maximal ideal
$\m \subset B$, and let $\wt{A}$ be a flat associative algebra over
$B$ such that $A = \wt{A}/\m$ is commutative. Then there exist a
canonical subalgebra $\wt{A}^p \subset \wt{A}$ such that $\wt{A}^p$
is flat over $B^p \subset B$, the natural map $\wt{A}^p
\otimes_{B^p} B \to \wt{A}$ is injective, and the quotient map
$\wt{A} \to A$ induces an isomorphism between $\wt{A}^p/\m^p \cong
(\wt{A}^p \otimes_{B^p} B)/\m$ and $A^p \subset A$.
\end{conj}

This conjecture does hold in the two extreme cases of commutative
deformations $\wt{A}$ ($\wt{A}^p$ is the subalgebra generated by
$p$-th powers) and of central quantizations of symplectic affine
manifolds ($\wt{A}^p$ is the center of the algebra
$\wt{A}$). Lemma~\ref{quasi.fr} is Conjecture~\ref{cnj} for $B =
k[h]/h^{p-1}$ (in this case $B^p = k$). Conjecture essentially says
that even in the general case, $A^p \subset A$ does extend to a flat
subalgebra in the whole $\wt{A}$, possibly no longer central nor
commutative, and this flat deformation of $A^p$ is moreover induced
from $B^p \subset B$.

If one thinks in terms of analogy between deformations and
automorphisms, then Conjecture~\ref{cnj} is analogous to the fact
that the automorphism group $\Aut A$ preserves $A^p \subset A$, and
its induced action on $A^p$ factors through the Frobenius map of the
natural group scheme structure. Of course, this analogy is flimsy at
best, so that we have no really convincing reason to believe that
the conjecture is true. If it is true, its proof might possibly
involve an extension of the restricted Lie algebra structure on
vector fields to some sort of structure on the whole Hochschild
cohomology DG Lie algebra of a commutative algebra $A$.

If one assumes Conjecture~\ref{cnj}, then the natural generalization
of Theorem~\ref{qua} would replace Frobenius-constant quantizations
with all quantizations $\wt{A}$ with prescribed algebra
$\wt{A}^p$. Alternatively, one might want to study quantizations
such that $\wt{A}^p \cong \wt{A}$. The latter would correspond to
the condition $C(\Omega)=\Omega$ in Theorem~\ref{cent}, and as we
noted in Remark~\ref{proj.rem}, it seems a more natural question for
compact manifolds -- such as even-dimensional abelian varieties, or
surfaces of type $K3$.

\section{Poisson structures.}

\subsection{Central quantizations.}\label{cent.sub}
The goal of this Section is to prove Theorem~\ref{cent} and
Theorem~\ref{obstr}. We start with some preliminaries on
differential geometry in positive characteristic.

Let $B$ be a Noetherian commutative algebra over the field $k$, and
let $A$ be a flat quasiregular commutative algebra over $B$.  Denote
by $A^{[p]} = A^p \otimes_{B^p} B \subset A$ the $B$-subalgebra
generated by $p$-th powers of elements in $A$; the algebra $A^{[p]}$
is a quotient of the pullback $A \otimes_{\Fr_B} B$ of the algebra
$A$ with respect to the Frobenius map $\Fr_B:B \to B$, and the
Frobenius map $\Fr_A:A \to A$ factors through a natural map map
$\Fr_{A^{[p]}}:A \to A^{[p]}$. Denote by $\Omega^\hdot(A) =
\Lambda^\hdot\Omega^1(A/B)$ the modules of differential forms on
$A/B$, and denote by $\Omega^\hdot_{cl}(A) \subset \Omega^\hdot(A)$
the subspaces of closed forms. Denote $\Omega^\hdot_{[p]}(A) =
\Omega^\hdot(A/B) \otimes_{\Fr_{A^{[p]}}} A^{[p]}$, and denote by
$C^\idot:\Omega^\hdot_{cl}(A) \to \Omega^\hdot_{[p]}(A)$ the
$A^{[p]}$-linear Cartier operations.

Consider the module $T(A) = T(A/B)$ of $B$-linear derivations of the
algebra $A$ (in geometric language, these are relative vector fields
on $\Spec A$ over $\Spec B$). It is a flat $A$-module and a Lie
algebra; moreover, for any derivation $\xi:A \to A$, its $p$-th
power $\xi^p:A \to A$ is also a derivation, so that $T(A)$ is a
restricted Lie algebra. The Lie algebra $T(A)$ acts on
$\Omega^\hdot(A)$ by the Lie derivative, and we will denote the
action of a vector field $\xi \in T(A)$ by $\ad
\xi:\Omega^\hdot(A)\to\Omega^\hdot(A)$.

We note that by extension of scalars every differential form $\alpha
\in \Omega^l(A)$ on $A/B$ induces a form $\Fr_{A^{[p]}}^*\alpha \in
\Omega^l_{[p]}(A)$. In particular, a function $f \in A$ is a
differential form of degree $0$, and we have $\Fr_{A^{[p]}}^*f = f^p
\in A^p$. Again by extension of scalars, every vector field $\xi \in
T(A)$ induces a contraction operation $\alpha \mapsto \alpha \cntrct
\xi$, $\Omega^{\hdot+1}_{[p]}(A) \to \Omega^\hdot_{[p]}(A)$.

We will need one fact on the relation between Cartier operations and
the restricted Lie algebra structure on $T(A)$. Namely, for any
vector field $\xi \in T(A)$ define an operation
$i_\xi^{[p]}:\Omega^{\hdot+1}(A) \to \Omega^\hdot(A)$ by
$$ 
i_\xi^{[p]}(\alpha) = \xi^{[p]} \cntrct \alpha - (\ad \xi)^{p-1}(\xi
\cntrct \alpha).
$$
By the Cartan homotopy formula, the operation $i_\xi^{[p]}$ commutes
with the de Rham differential, so that it induces an operation on
the de Rham cohomology groups $H^\hdot_{DR}(A) \cong
\Omega^\hdot(A^p)$.

\begin{lemma}\label{car.p.le}
For every closed form $\alpha \in \Omega^\hdot_{cl}(A)$ we have
\begin{equation}\label{car.p}
C^{\hdot-1}\left(i_\xi^{[p]}(\alpha)\right) = C^\idot(\alpha) \cntrct \xi.
\end{equation}
\end{lemma}

\proof{} Both sides are additive and $A^p$-linear with respect to
$\alpha$. Moreover, for any closed forms $\alpha,\beta$ we have
$$
(\xi \cntrct \alpha) \wedge \ad \xi(\beta) = (-1)^{\deg
\alpha}\ad \xi(\alpha) \wedge (\xi \cntrct \beta) - d((\xi \cntrct
\alpha) \wedge (\xi \cntrct \beta);
$$
therefore
$$
\begin{aligned}
i_\xi^{[p]}(\alpha \wedge \beta) &- i_\xi^p\alpha \wedge \beta + \alpha
\wedge i_\xi^p\beta =\\
&=\sum_{i=1}^{p-2}\binom{p-1}{i}(\xi \cntrct (\ad \xi)^i\alpha)
\wedge (\ad \xi)^{p-1-i}\beta \\
&\qquad + (-1)^{\deg
\alpha}\sum_{i=1}^{p-2}\binom{p-1}{i}(\ad \xi)^i\alpha \wedge (\xi
\cntrct (\ad \xi)^{p-1-i}\beta)\\
&= \sum_{i=1}^{p-2}\left(\binom{p-1}{i} +
\binom{p-1}{i+1}\right)(\xi \cntrct (\ad \xi)^i\alpha) \wedge (\ad
\xi)^{p-1-i}\beta \\
&\qquad + d((\ad \xi)^{p-2}((\xi \cntrct \alpha) \wedge
(\xi \cntrct \beta))\\
&=d((\ad \xi)^{p-2}((\xi \cntrct \alpha) \wedge (\xi \cntrct
\beta)),
\end{aligned}
$$
which implies that the map $H^{\hdot+1}_{DR}(A) \to H^\hdot_{DR}(A)$
induced by $i_\xi^{[p]}$ is a derivation with respect to the algebra
structure on $H^\hdot_{DR}$. Thus it suffices to prove \eqref{car.p}
for forms $\alpha$ of degree $1$. Let $\alpha$ be such a form,
consider the connection $\nabla = d + \alpha$ on the $A$-module $A$,
and denote by $f:A \to A$ multiplication by $\alpha \cntrct \xi$. By
definition of the Cartier operator, $C^1(\alpha) -
\Fr_{A^{[p]}}^*\alpha$ is the $p$-curvature of the connection
$\nabla$; in other words, we have
$$
(\nabla_\xi)^p = \nabla_{\xi^{[p]}} + C^1(\alpha) \cntrct \xi -
\Fr_{A^{[p]}}^*(\alpha \cntrct \xi) = \nabla_{\xi^{[p]}} + C^1(\alpha)
\cntrct \xi - f^p.
$$
Since $\nabla_\xi = \xi + f$, \eqref{car.p} is therefore equivalent to
$$
(\xi+f)^p = \xi^p + \alpha \cntrct \xi^{[p]} + f^p -
i_\xi^{[p]}(\alpha).
$$
This immediately follows from Lemma~\ref{sq}: take $x=\xi$, $y=f$,
and consider the Lie algebra $L$ of differential operators $A \to A$
of order $\leq 1$ with the natural filtration by order (inverted and
shifted by $1$).
\endproof

Assume now that the quasiregular algebra $A/B$ is equipped with a
non-degenerate closed $2$-form $\Omega \in \Omega^2_{cl}(A)$. As in
characteristic $0$, the form $\Omega$ induces a $B$-linear Poisson
bracket on the algebra $A$. This structure is non-degenerate in the
sense that the Poisson bivector $\Theta \in \Lambda^2(T(A))$ is a
non-degenerate bilinear form on the $A$-module $\Omega^1(A)$;
conversely, every non-degenerate $B$-linear Poisson bracket on $A$
comes from a symplectic form. For any $f \in A$ we denote the
derivation $a \mapsto \{a,f\}$ by $H_f$; explicitly, we have $H_f
\cntrct \Omega = da$. Vector fields $H_f$ are called Hamiltonian and
form a Lie subalgebra $H(A) \subset T(A)$ in the Lie algebra $T(A)$.

\begin{corr}\label{car.0}
The subalgebra $H(A)$ is closed under the restricted Lie algebra
structure if and only if $C^2(\Omega) = 0$.
\end{corr}

\proof{} By definition, for any $f \in A$ we have $C^1(H_f \cntrct
\Omega) = C^1(df) = 0$. Lemma~\ref{car.p.le} yields
$$
C^1(H_f^{[p]} \cntrct \Omega) = H_f^\tw \cntrct C^2(\Omega).
$$
The vector field $H_f^{[p]}$ is Hamiltonian if and only if the
left-hand side vanishes; therefore $H(A) \subset T(A)$ is closed
under the restricted Lie algebra structure if and only if the
right-hand side vanishes for any $f$. Since $\Omega$ is
non-degenerate, Hamiltonian vector fields generate $T(A)$ as an
$A$-module, so that this is possible if and only if $C^2(\Omega)=0$.
\endproof

\proof[The proof of Theorem~\ref{cent}.]\label{cent.pf} All claims
are local, therefore we can assume that $X = \Spec A$ and $S = \Spec
B$ are affine. Then equivalence between \thetag{i} and \thetag{ii}
becomes Corollary~\ref{car.0}. To prove the second claim, we can
assume that we are given a quasiregular $B$-algebra $A$ equipped with
a non-degenerate closed $2$-form $\Omega \in \Omega^2_{cl}(A)$, and
its quantization $A_h$ such that the center $\ZZ_h \subset A_h$
satisfies $\ZZ_h/h \cong \ZZ \subset A$. Moreover, by
Lemma~\ref{no.poi.cent} the Poisson center $\ZZ \subset A$ coincides
with the subalgebra $A^{[p]} \subset A$.  Consider the quotient
$A_h/h^p$, take an element $f \in A$ and lift $f^p \in A^{[p]}$ to
an element $\wt{f}^p \in \ZZ_h$. By Lemma~\ref{quasi.fr} we have
$$
\wt{f}^p = f^p + f_1h + \dots + f_{p-2}h^{p-2} + gh^{p-1} \mod
h^p
$$
for some $f_1,\dots,f_{p-2} \in A^{[p]}$ and $g \in A$. By
\eqref{ad.fr}, $[f_i,y]$ is divisible by $h^{p-1}$ for all $f_i$ and
any $y \in A$, and the same is true for $f^p$. On the other hand,
since $\wt{f}^p$ is central, we have $[\wt{f}^p,y]=0$; reducing this
modulo $h^{p+1}$, we obtain
$$
(\ad f)^p(y)+[g,y]h^{p-1} = 0 \mod h^{p+1},
$$
which further reduces to $H_f^{[p]}(y) = H_g(y)$ and $H_f^{[p]} =
H_g$. This finishes the proof.
\endproof

\subsection{Conformal derivations.}
Having proved Theorem~\ref{cent}, we now proceed to
Theorem~\ref{obstr}.

\begin{defn}
Assume given a Poisson algebra $A/k$ and a constant $\lambda \in
k$. A {\em conformal derivation} $\xi$ of the algebra $A$ with {\em
weight $\lambda$} is a $k$-linear map $\xi:A \to A$ which is a
derivation with respect to the multiplication in $A$ and satisfies
\begin{equation}\label{conf.poi}
\xi(\{x,y\}) = \{\xi(x),y\} + \{x,\xi(y)\} - \lambda\{x,y\}.
\end{equation}
\end{defn}

A conformal derivation of weight $0$ is just a Poisson derivation,
in particular, all Hamiltonian derivations are conformal of weight
$0$.

\begin{lemma}
Assume given a Poisson algebra $A/k$ and a conformal derivation
$\xi:A \to A$ of some weight $\lambda \in k$. Then for any two
elements $x,y \in A$ we have
\begin{equation}\label{conf.poly}
\begin{aligned}
\xi(L(x,y)) &= \lambda L(x,y) + \ad(x+y)^{p-1}(\xi(x+y))\\
&\quad -
\ad(x)^{p-1}(\xi(x)) - \ad(y)^{p-1}(\xi(y)),\\
\xi(P(x,y)) &= \lambda P(x,y) + \ad(xy)^{p-1}(\xi(xy))\\
&\quad -
x^p\ad(y)^{p-1}(\xi(y)) - y^p\ad(x)^{p-1}(\xi(x)),
\end{aligned}
\end{equation}
where we denote $\ad(z)(w) = \{z,w\}$.
\end{lemma}

\proof{} We pass to the universal situation. Consider the algebra
$Q^\hdot(x,y,z)$ of quantized polynomials in three variables
$x,y,z$. By definition, the localization algebra
$Q^\hdot(x,y,z)(h^{-1})$ is just the free associative
$k[h,h^{-1}]$-algebra generated by $x$, $y$ and $z$; therefore there
exists a unique derivation $e:Q^\hdot(x,y,z)(h^{-1}) \to
Q^\hdot(x,y,z)(h^{-1})$ which annihilates $x$, $y$ and $z$ and
satisfies $e(h)=\lambda h$. This derivation obviously preserves the
subalgebra $Q^\hdot(x,y,z) \subset Q^\hdot(x,y,z)(h^{-1})$. Consider
the derivation $\xi = e+ \ad(z)$ of the algebra
$Q^\hdot(x,y,z)$. Define a grading on $Q^\hdot(x,y,z)$ by setting
$\deg x = \deg y = \deg h = 0$, $\deg z = 1$, and let
$\wt{Q}^\hdot(x,y)$ its the component of degree $0$. The derivation
$\xi$ and the quantized algebra structure are compatible with the
grading and restrict to $\wt{Q}^\hdot(x,y) \subset Q^\hdot(x,y,z)$.

The quotient $\overline{Q}^\hdot(x,y) = \wt{Q}^\hdot(x,y)/h$ of the
algebra $\wt{Q}^\hdot(x,y)$ is a Poisson algebra equipped with a
derivation $\xi$, and since $e(h)=\lambda h$, the derivation $\xi$
is conformal of weight $\lambda$. It is easy to check that the
Poisson algebra $\overline{Q}^\hdot(x,y)$ is freely generated by
$\xi^l(x)$ and $\xi^l(y)$, $l \geq 0$. Therefore for any two
elements $x,y \in A$ in a Poisson algebra $A$ equipped with a
conformal derivation $\xi$ of weight $\lambda$, there exists a
unique Poisson map $s:\overline{Q}^\hdot(x,y) \to A$ which sends $x$
to $x$, $y$ to $y$, and satisfies $w(\xi(a))=\xi(w(a))$ for any $a
\in \overline{Q}^\hdot(x,y)$.

We conclude that it suffices to prove the claim for $A =
\overline{Q}^\hdot(x,y)$. Moreover, it even suffices to prove it for
the quantized subalgebra $\wt{Q}^\hdot(x,y) \subset
Q^\hdot(x,y,z)$. But since the algebra $Q^\hdot(x,y,z)$ has no
$h$-torsion, we may multiply both sides of \eqref{conf.poly} by
$h^{p-1}$; both equations then immediately follow from \eqref{ad.fr}
and \eqref{poly}.
\endproof

\subsection{Restricted Poisson structures.}\label{obstr.sub}
Return to the notation of Subsection~\ref{cent.sub}: $A$ is a
quasiregular algebra over $B$ with a non-degenerate $B$-linear
Poisson structure induced by a symplectic form $\Omega \in
\Omega^2(A)$. By Lemma~\ref{no.poi.cent}, the subalgebra $A^{[p]}
\subset A$ is the Poisson center of the algebra $A$.

\begin{lemma}\label{ravno.lemma}
Assume given a restricted structure on the Poisson algebra $A$. Then
a Poisson derivation $\xi \in T(A)$ is Hamiltonian if and only if we
have
\begin{equation}\label{ravno}
\xi(a^{[p]}) = \ad(a)^{p-1}(\xi(a))
\end{equation}
for any $a \in A$.
\end{lemma}

\proof{} Since $\xi$ is Poisson, the $1$-form $\beta = \xi \cntrct
\Omega$ is closed, $C^1(\beta)$ is well-defined, and $\xi$ is
Hamiltonian if and only if $C^1(\beta)=0$. Since Hamiltonian vector
fields generate the $A$-module $T(A) = T(A/B)$, this holds if and
only if
$$
H_a^\tw \cntrct C^1(\beta) = 0
$$
for any Hamiltonian vector field $H_a$, $a \in A$. By
Lemma~\ref{car.p.le}, this is equivalent to
$$
C^0\left(i_{H_a}^{[p]}(\beta)\right) = 0,
$$
and since the Cartier operation $C^0$ is injective, this reduces to
$i_{H_a}^{[p]}(\beta)=0$, which by definition exactly means
\eqref{ravno}. 
\endproof

\begin{prop}\label{conf.prop}
Assume given a conformal derivation $\xi:A \to A$ of weight
$\lambda$ of a Poisson algebra $A$.
\begin{enumerate}
\item For any restricted structure on $A$, the map $\kappa:A \to A$
given by
\begin{equation}\label{kappa=>p}
\kappa(a) = (\xi - \lambda\id)(a^{[p]}) - \ad(a)^{p-1}(\xi(a))
\end{equation}
is a Frobenius-derivation of the algebra $A$ with values in the
Poisson center $A^{[p]} \subset A$.
\item If $\lambda=1$, then every Frobenius-derivation $\kappa:A \to
A^{[p]} \subset A$ is given by \eqref{kappa=>p} for a unique
restricted structure on the Poisson algebra $A$.
\end{enumerate}
\end{prop}

\proof{} To prove \thetag{i}, we apply \eqref{conf.poly} and deduce
that $\kappa(a+b)=\kappa(a)+\kappa(b)$ and $\kappa(ab) =
a^p\kappa(b) + b^p\kappa(a)$. Moreover, for any $a,b \in A$ we have
$$
\begin{aligned}
\{\kappa(a),b\} &= \xi(\{a^{[p]},b\}) - \{a^{[p]},\xi(b)\} -
	[\xi,H_a^{[p]}](b)\\
                &= \xi(H_a^{[p]}(b))-H_a^{[p]}(\xi(b))
	-[\xi,H_a^{[p]}](b) =0,
\end{aligned}
$$
so that $\kappa(a)$ indeed lies in the Poisson center $A^{[p]}
\subset A$.
                
To prove \thetag{ii}, denote $\alpha = \xi \cntrct \Omega$. By the
Cartan Homotopy formula, $d\alpha=\lambda\Omega$, and
\eqref{kappa=>p} can be rewritten as
$$
a^{[p]} = i_{H_a}^{[p]}\alpha - \kappa(a).
$$
Therefore the map $a \mapsto a^{[p]}$ is indeed uniquely defined by
$\xi$ and $\kappa$, and we have to prove that it satisfies
\eqref{restr.lie.1}, \eqref{restr.lie.2}, and \eqref{restr.poi}. The
equation \eqref{restr.lie.1} is immediate: since $d\kappa(a)=0$, we
have $da^{[p]}= di_{H_a}^{[p]}\alpha = H_a^{[p]} \cntrct \Omega$. To
prove the rest, note that any function $f \in A$ which commutes with
all Hamiltonian derivations lies in the Poisson center $A^{[p]}
\subset A$. For such a function, we have $\xi(f)=0$ and
$(\xi-\id)(f)=-f$. Therefore it suffices to prove both
\eqref{restr.lie.2} and \eqref{restr.poi} after applying to both
sides the map $(\xi-\id)$ and the maps $H_a$ for all $a \in A$. Both
equations then immediately follow from \eqref{kappa=>p} and
\eqref{conf.poly}.
\endproof

\proof[Proof of Theorem~\ref{obstr}.]\label{obstr.pf} The claim is
local in Zariski topology, so that we can assume that $S=\Spec B$
and $X = \Spec A$ are affine, with $A$ quasiregular over $B$. We are
given a Frobenius-derivation $K:B \to B$, and we have to establish a
bijection between the set of all restricted structures on the
Poisson algebra $A$ which are compatible with $K:B \to B$ and the
set of all isomorphisms between $\Omega^1_{[p]}(A)$-torsors
$C^1_*[\M_\Omega]$ and $\Fr^*_{A^{[p]}}\Omega \cntrct
[\theta_K]$. By the Cartan Homotopy formula, sections of the torsor
$[\M_\Omega]$ correspond to conformal derivations $\xi:A \to A$ of
weight $1$. By Proposition~\ref{conf.prop}~\thetag{i}, every
restricted structure on $A$ induces a map from the set of all
conformal derivations $\xi:A \to A$ of weight $1$ to the set of all
Frobenius-derivations $\kappa:A \to A^{[p]} \subset A$. It
immediately follows from \eqref{kappa=>p} that no matter what $\xi$
one takes, a restricted structure is compatible with $K:B \to B$ if
and only if the corresponding $\kappa:A \to A^{[p]} \subset A$ is
equal to $K$ on $B \subset A$.

But the $B$-modules of $B$-linear Frobenius derivations $A \to
A^{[p]} \subset A$ is naturally identified with $\Omega^1_{[p]}(A)$
by $\xi \mapsto \Fr^*_{A^{[p]}}\Omega \cntrct \xi$. Therefore the
torsor $\Fr^*_{A^{[p]}}\Omega \cntrct [\theta_K]$ is naturally
identifies with the set of all Frobenius-derivations $\kappa:A \to
A^{[p]} \subset A$ which are equal to $K$ on $B \subset A$.

In other words, a restricted structure is compatible with $K$ if and
only if the correspondence $\xi \mapsto \kappa$ actually defines a
map $[\M_\Omega] \to \Fr^*_{A^{[p]}}\Omega^\tw \cntrct
[\theta_K]$. By Lemma~\ref{ravno.lemma}, two conformal derivations
give the same $\kappa$ if and only if their difference is a
Hamiltonian derivation; by the Cartan Homotopy Formula, this means
that this map actually factors through the quotient
$C^1_*[\M_\Omega] = [\M_\Omega]/d(A)$ and gives an isomorphism
$C^1_*[\M_\Omega] \cong \Omega^\tw \cntrct [\theta_K]$. Finally, by
Proposition~\ref{conf.prop}~\thetag{ii}, every such isomorphism
comes from a unique restricted structure on $A$ compatible with $K:B
\to B$.
\endproof

\begin{remark}
Informally, the situation can be described as follows. There are
three natural $\Omega^1_{[p]}(A)$-torsors associated to $A/B$ and
$K:B \to B$. Firstly, there is a torsor of all $B$-linear conformal
derivations of weight $1$, considered modulo Hamiltonian
derivations. This torsor does not depend on $K$. Secondly, there is
the torsor $\Fr^*_{A^{[p]}}\Omega \cntrct [\theta_K]$, which is
identified with the set of all Frobenius-derivations $A \to A^{[p]}
\subset A$ which are equal to $K$ on $B \subset A$. This torsor, in
spite of notation, does not depend on the symplectic form
$\Omega$. Finally, there is the set of restricted structures on $A$
compatible with $K:B \to B$, and this set is also a torsor over
$\Omega^1_{A^{[p]}}$ -- recall (page~\pageref{dff}, remarks after
Definition~\ref{restr.qua.def}) that the difference of two
restricted structures is a Frobenius-derivation $A \to A$ with
values in the Poisson center $A^{[p]} \subset
A$. Proposition~\ref{conf.prop} shows that the sum of these three
torsors is canonically trivialized. Therefore giving an element in
one of them is equivalent to giving an isomorphism between the other
two.
\end{remark}

\section{Local theory.}\label{loc.sec}

To proceed further, we need to study symplectic manifolds and their
quantizations locally. In characteristic $0$, ``locally'' means at
best ``in the formal neighborhood of a point''. In positive
characteristic, we can do better and replace a regular complete
local algebra $A'$ with its quotient $A= A'/\m^{[p]} \cdot A'$ by
the $p$-th powers of elements in the maximal ideal $\m \subset
A'$. The algebra $A$ is a finite-dimensional local Artin algebra
over $k$, and it is quasiregular. The goal of this section is to
prove Proposition~\ref{step} and Proposition~\ref{level.1} -- or
rather, more precise versions of them which take account of the
category structure -- for $X = \Spec A$.

\subsection{Differential geometry.}\label{drb.sub}

Fix a finite-dimensional vector space $W$ over $k$. Consider the
algebra $A = k[[W^*]]/\m^{[p]}$, where $\m^{[p]} \subset k[[W^*]]$
is the ideal generated by $a^p$ for all $a$ in the augmentation
ideal $\m \subset k[[W^*]]$. Then $A$ is a finite-dimensional local
algebra over $k$, of rank $p^{\dim W}$, with maximal ideal $\m
\subset A$. Moreover, the algebra $A$ is quasiregular: the module
$T(A) \cong W \otimes A$ of derivations of the algebra $A$ is a free
$A$-module generated by the vector space $W$. Analogously, for every
$k \geq 0$, the module $\Omega^k(A) \cong \Lambda^k(W^*) \otimes A$
of $k$-forms on $A$ is a free $A$-module generated by
$\Lambda^k(W^*)$. The algebra $T(A)$ is a restricted Lie algebra;
the Cartier operations $C^k$ induce graded ring isomorphisms
$H^k_{DR}(A) \cong \Omega^k_{[p]}(A) \cong \Lambda^k(W^{*\tw})$.

The group $\Aut A$ of all automorphisms of the algebra $A$ has a
natural structure of a quasiregular group scheme over $k$, with Lie
algebra $T(A)$. The subgroup $(\Aut A)_0 \subset \Aut A$ of all
automorphisms that preserve $\m \subset A$ is a regular group scheme
with Lie algebra $T(A)_0 = W \otimes \m \subset W \otimes A \cong
T(A)$. The quotient homogeneous space $(\Aut A)/(\Aut A)_0$ is
naturally identified with $\Spec A$.

Assume that $A$ is equipped with a non-degenerate closed $2$-form
$\Omega \in \Omega^2(A)$ (in particular, $W$ is
even-dimensional). The form $\Omega$ defines a Poisson structure on
the algebra $A$ and the subalgebra $H(A) \subset T(A)$ of
Hamiltonian vector fields. We record the following obvious facts.

\begin{lemma}\label{poi.simple}
The only Poisson ideals in $A$ are $0$ and $A$ itself.
\end{lemma}

\proof{} Since $A$ is local, any proper ideal $I \subset A$ lies in
$\m^k$ for some $k$; but since Hamiltonian vector fields generate
$T(A)$ as an $A$-module, the ideal $I$ must be preserved by all
derivations of the algebra $A$. An easy induction on $k$ shows that
this is possible only for $I=0$.
\endproof

\begin{lemma}\label{loc.azu}
For any quantization $A_h$ of the algebra $A$, the algebra
$A_h(h^{-1})$ is an Azumaya algebra over the Laurent power series
field $k((h))$.
\end{lemma}

\proof{} By Lemma~\ref{no.poi.cent}, the Poisson center of the
algebra $A$ is $k \subset A$; therefore the algebra $A_h(h^{-1})$ is
central, and it remains to prove that it is a simple
algebra. Indeed, let $I \subset A_h(h^{-1})$ be a two-sided
ideal. Then $I_0 = I \cap A_h$ is a two-sided ideal in $A_h \subset
A_h(h^{-1})$, and $I = I_0(h^{-1})$. We further note that by
definition, $I_0/h \subset A$ is a Poisson ideal in $A$. By
Lemma~\ref{poi.simple}, this means that either $I_0/h=0$ or
$I_0/h=A$. By the Nakayama Lemma, this implies that $I_0=0$ or
$I_0=A_h$, which in turn yields $I=0$ or $I=A_h(h^{-1})$.
\endproof

Assume now that we are given a local Artin $k$-algebra $B$ equipped
with a Frobenius-derivation $K:B \to B$, and assume that the Poisson
algebra $A \otimes B$ over $B$ is equipped with a restricted
structure compatible with $K$. Denote by $\Aut(A)^B$ the algebraic
group of all $B$-linear automorphisms of the algebra $A \otimes B$;
it is obtained from $\Aut(A)$ by the Weil restriction of scalars
from $B$ to $k$. The subgroup $\G^B \subset \Aut(A \otimes B/B)$ of
automorphisms that preserve the (relative) symplectic form and the
restricted structure inherits a structure of an algebraic group. We
also have the subgroup $\Aut(A)^B_0 \subset \Aut(A)^B$ of
automorphisms that preserve the maximal ideal in $A \otimes B$;
denote the intersection $\Aut(A)^B_0 \cap \G^B \subset \Aut(A)^B$ by
$\G_B^0$.

\begin{lemma}\label{ham}
A $B$-linear derivation $\xi:A \otimes B \to A \otimes B$ lies in
the Lie algebra of the subgroup $\G^B$ if and only if it is
Hamiltonian, $\xi = H_f$ for some $f \in A \otimes B$.
\end{lemma}

\proof{} To compute the Lie algebra of $\G^B$, consider the dual
numbers algebra $k \langle \eps \rangle = k[\eps]/\eps^2$. By
definition, a derivation $\xi$ comes from the Lie algebra of $\G^B$
if and only if $1+\xi\eps:A \otimes B\langle \eps \rangle \to A
\otimes B\langle \eps \rangle$ is compatible with the symplectic
form and the restricted structure. Compatibility with the symplectic
form means that $\xi$ is a Poisson derivation. Compatibility with
the restricted structure reads as
\begin{align*}
a^{[p]} + \xi(a^{[p]})\eps &= (a + \xi(a)\eps)^{[p]}\\
&= a^{[p]} + L(a,\xi(a)\eps) + \xi(a)^{[p]}\eps^p + \xi(a)^p\eps^{[p]} +
P(\xi(a),\eps)
\end{align*}
for any $a \in A \otimes B$. By definition, $\eps^p=\eps^{[p]}=0$;
moreover, since $\{-,\eps\}$ vanishes tautologically, we have
$P(\xi(a),\eps)=0$. Therefore only the first two terms in the
right-hand side are possibly non-trivial. Applying Lemma~\ref{sq} to
$A \otimes B \cdot \eps \subset A \otimes B\langle\eps\rangle$, we
see that this is equivalent to \eqref{ravno}, which by
Lemma~\ref{ravno.lemma} holds if and only if $\xi$ is Hamiltonian.
\endproof

\begin{prop}\label{darboux}
Assume given an Artin $k$-algebra $B$ and a Frobenius-derivation
$K:B \to B$. Then all the non-degenerate restricted Poisson
structures on the algebra $A \otimes B$ over $B$ compatible with $K$
are isomorphic. Moreover, the group $\G^B_0 \subset \Aut(A)^B$ is a
regular algebraic group.
\end{prop}

\proof{} The set $R$ of all non-degenerate restricted Poisson
structures on $A \otimes B$ over $B$ has a natural structure of an
algebraic variety, and this variety is irreducible. Indeed,
forgetting the restricted structure defines a map $R \to U \subset
\Omega^2(A \otimes B)_{ex}$ onto the set $U$ of all non-degenerate
$B$-linear Poisson structures on $A \otimes B$ which admit a
restricted power operation, and by Theorem~\ref{cent} the set $U$ is
an open subset in the affine space of exact relative $2$-forms on $A
\otimes B$ over $B$; on the other hand, by Theorem~\ref{obstr} $R$
is an affine torsor over $U$ (with respect to $\Omega^1_{[p]}(A
\otimes B)$ considered as an additive algebraic group). Now, the
group $\Aut(A)_0^B$ acts on $R$, and we have to show that this
action is transitive. Since $U$ is irreducible, it suffices to show
that the orbit $\Aut(A)_0^B \cdot r$ of any point $r \in R$ is the
whole $R$. But by definition the stabilizer of a point $r \in R$ is
$\G_0^B$, the automorphism group of $A \otimes B$ with the
corresponding restricted Poisson structure. Therefore
\begin{equation}\label{dimm}
\begin{aligned}
\dim R &\geq \dim \Aut(A)^B_0 \cdot r = \dim \Aut(A)^B_0 - \dim
\G^B_0\\ 
&\geq \dim \Aut(A)^B_0 - \dim H_0(A \otimes B),
\end{aligned}
\end{equation}
where $H_0(A \otimes B)$ is the Lie algebra of the group $\G^B_0$ --
which by Lemma~\ref{ham} coincides with the Lie algebra of
Hamiltonian vector fields on $A \otimes B$ preserving the
maximal ideal. But $\dim R = \dim U + \dim \Omega^1_{[p]}(A \otimes
B) = \dim \Omega^2_{ex}(A \otimes B) + \dim H^1(\Omega^\hdot(A
\otimes B)) = \dim \Omega^1(A \otimes B) - \dim \Omega^1_{ex}(A
\otimes B)$, while on the other hand, $\dim \Aut(A)^B_0 - \dim H_0(A
\otimes B) = \dim T(A \otimes B) - \dim H(A \otimes B) = \dim
\Omega^1(A \otimes B) - \dim \Omega^1_{ex}(A \otimes B)$. Therefore
all the inequalities in \eqref{dimm} are in fact equalities. This
means, firstly, that $\dim R = \dim \Aut(A)_0^B \cdot r$, so the
both sets coincide, and secondly, that $\dim \G^B_0$ is the same as
the dimension of its Lie algebra, which means that the group
$\G^0_B$ is regular.
\endproof

\begin{remark}
If $B=k$, this Proposition in particular implies that all symplectic
forms $\Omega \in \Omega^2(A)$ with $C^2(\Omega)=0$ are isomorphic
(every such form admits a compatible restricted power operation by
Theorem~\ref{obstr}). Since $C^2(\Omega)$ is an obvious invariant of
a symplectic form, this is the best replacement for the Darboux
Theorem one could hope for in characteristic $> 0$.
\end{remark}

\subsection{Quantizations.}\label{qua.sub}

We can now prove Proposition~\ref{step} for $X = \Spec A$. Assume
given a quantization base $B_0$, a small extension $I \subset B \to
B_0$, and a $B_0$-quantization $A_0$ of the algebra $A$. Denote by
$\m_B \subset B$ the maximal ideal of the local $k$-algebra
$B$. Consider the category $Q(A_0,B)$ of all regular
$B$-quantizations $A_B$ of the algebra $A$ equipped with an
isomorphism $A_B/IA_B \cong A_0$. (We remind the reader that a
$B$-quantization $A_B$ is called regular if and only if any $a \in
A_B$ such that $a^p=0 \mod \m_BA_B$ satisfies $a^p =
h^{p-1}a^{[p]}$.)

\begin{prop}\label{small.ex}
All objects in the category $Q(A_0,B)$ are isomorphic.
\end{prop}

\proof{} Assume first that $B_0=k$, so that $I \subset B$ coincides
with the maximal ideal $\m_B \subset B$ and satisfies
$\m_B^2=0$. Assume also that $h=0$, so that a quantization $A_B
\subset Q(A_0,B)$ is a restricted Poisson algebra and a square-zero
extension of the restricted Poisson algebra $A_0=A$. We have to show
that $A_B \cong A \otimes_k B$.  Indeed, the algebra $A$ is by
definition generated by $k=\dim W$ elements $a_1,\dots,a_k \in A$
modulo the relations $a_i^p=0$. Lift the elements $a_i$ to elements
$\wt{a}_i \in A_B$; since the quantization $A_B$ is regular, we
still have $\wt{a}_i^p=h^{p-1}\wt{a}_i^{[p]}=0$. This means that
$A_B \cong A \otimes B$ as commutative $B$-algebras. But by
Proposition~\ref{darboux} all non-degenerate restricted Poisson
algebra structures on $A \otimes B$ are isomorphic. It remains to
notice that since $A_B/\m_B$ is isomorphic to the non-degenerate
Poisson algebra $A$, the Poisson structure on $A_B$ is also
non-degenerate.

Consider now the general case. For any two algebras $C_1$, $C_2$
equipped with surjective maps $C_1 \to C$, $C_2 \to C$ to a third
algebra $C$, denote by $C_1 \oplus^C C_2$ their Baer sum -- that is,
the subalgebra in $C_1 \oplus C_2$ of elements with the same image
in $C$. Assume given two different regular $B$-quantizations
$A_1,A_2 \in Q(A_0,B)$. Let $B'=B \oplus^{B_0} B$, and let $A_{12} =
A_1 \oplus^{A_0} A_2$. The algebra $A_{12}$ is obviously a regular
$B'$-quantization of the algebra $A$. Let $p_1:B'=B \oplus^{B_0} B
\to B$ be the projection onto the first summand, let $\delta:B \to
B'$ be the diagonal embedding, and let $\tau:B' \to B''$ be the
natural projection onto the quotient $B'' = B'/\delta(\m_B)$. Then
we obviously have $B'' \cong k \oplus I$, and the map $p_1 \oplus
\tau:B' \to B \oplus B''$ identifies the algebra $B'$ with the Baer
sum $B \oplus^{k} B''$. Under this identification, $A_{12}/(0 \oplus
I)A_{12}$ is isomorphic to $A_1$, $A_{12}'=A_{12}/(\m_B \oplus
0)A_{12}$ is a regular $B''$-quantization of the algebra $A$, and we
have $A_{12} = A_1 \oplus^{A_0} A_{12}'$.

Now, we have $I^2=0$ and $h=0$ in the algebra $B'' =
B'/\delta(\m_B)$; therefore, as we have proved already, $A_{12}'$
must be isomorphic to $A \otimes B'' \cong A \oplus IA$ as a
restricted Poisson algebra. Therefore $A_{12} \cong A_1 \oplus IA$
as restricted Poisson algebras. In particular, we have a restricted
Poisson map $A_1 \to A_{12}$; composing it with the natural
projection $A_{12} \to A_2$ gives the required isomorphism $A_1
\cong A_2$.  \endproof

\begin{corr}\label{iso.lifts}
Assume given $A_1 \in Q(A_0,B)$, and consider the algebraic group
$\Aut(A_1)$ of all $B$-linear automorphisms of the restricted
quantized algebra $A_1$. Then the natural reduction map $\Aut(A_1)
\to \Aut(A_0)$ is a surjective map of algebraic groups over $k$.
\end{corr}

\proof{} It suffices to prove that the map $\Aut(A_1)(B') \to
\Aut(A_0)(B')$ of groups of $B'$-points is surjective for any
$k$-algebra $B'$. It is well-known that it suffices to prove it when
$B'$ is an arbitrary local Artin algebra, with maximal ideal
$\m_{B'} \subset B'$. By definition we have $\Aut(A_1)(B') \cong
\Aut(A_1 \otimes B')$ and $\Aut(A_0)(B') \cong \Aut(A_0 \otimes
B')$. Filtering $I \otimes B' \subset B \otimes B'$ by ideals $I
\otimes \m_{B'}^k \subset I \otimes B'$, we decompose the extension
$B \otimes B' \to B_0 \otimes B'$ into a series of small extensions.
By induction, it suffices therefore to prove that for any small
extension, the reduction map is surjective on the groups of points.
In other words, we have to prove that any $g \in \Aut(A_0)$ can be
lifted to an element $\wt{g} \in \Aut(A_1)$.

Indeed, let $A_2 \in Q(A_0,B)$ be the algebra $A_1$ equipped with a
different map $A_1 \to A_0$ -- namely, the composition of the
original one with the automorphism $g:A_0 \to A_0$. By
Proposition~\ref{small.ex}, $A_1$ and $A_2$ are isomorphic in
$Q(A_0,B)$; the isomorphism is the desired lifting $\wt{g} \in
\Aut(A_1)$.
\endproof

As an application of Proposition~\ref{small.ex}, we can extend
Proposition~\ref{darboux} to the following uniqueness result.

\begin{prop}\label{q.uni}
Assume given a quantization base $B$. There exists a unique, up to
an isomorphism, regular $B$-quantization $A_B$ of the Poisson
algebra $A$.
\end{prop}

\proof{} Uniqueness follows immediately by induction from
Proposition~\ref{small.ex}. Moreover, it suffices to prove existence
for $B = k[[h]]$: for an arbitrary $B$, we can take $A_B =
A_{k[[h]]} \otimes_{k[[h]]} B$. But by Proposition~\ref{darboux} we
can assume that for some integer $n$, $A$ is the commutative algebra
over $k$ generated by $x_1,\dots,x_n,y_n,\dots,y_n$ modulo the
relations $x_i^p=y_j^p=0$, and the symplectic form $\Omega$ on $A$
is given by
$$
\Omega = dx_1 \wedge dy_1 + \dots + dx_n \wedge dy_n.
$$
Then a regular $k[[h]]$-quantization of the algebra $A$ is given
by the so-called {\em reduced Weyl algebra} $D$, which is generated
by $x_1,\dots,x_n,y_1,\dots,y_n$ modulo the relations
$$
\begin{aligned}
x_ix_j - x_jx_i &= y_iy_j - y_jy_i = 0,\\
x_iy_j - y_jx_i &= \delta_{ij}h,\\
x_i^p = y_j^p &= 0,
\end{aligned}
$$
where $\delta_{ij}$ is the Kronecker delta symbol. It is immediately
obvious that this quantization is regular and Frobenius-constant
(the splitting map $s$ vanishes on all the generators $x_i$, $y_j$).
\endproof

\subsection{Automorphism groups.}\label{aut.sub}

We finish this section with a more detailed description of the
automorphism groups of quantizations. By Proposition~\ref{q.uni},
for any quantization base $B$ there exists a unique regular
$B$-quantization of the algebra $A$; denote this $B$-quantization by
$A_B$. Denote also by $\G^B = \Aut(A_B)$ the algebraic group of all
$B$-linear automorphisms of the $B$-quantization $A_B$. If $B=k$, we
have $\G^k=\HH$, the group of Hamiltonian automorphisms of the
restricted Poisson algebra $A$.

As in Subsection~\ref{exi.sub}, for any small Dieudonn\'e module $I$
over $k$ we denote by $\hh(A,I)$ the kernel of the map $C^1 -
\Fr_{A^{[p]}}^* \otimes K_I:\Omega^1_{cl}(A) \otimes I \to
\Omega^1_{[p]}(A) \otimes I$, so that we have an exact sequence
\begin{equation}\label{seq.3}
\begin{CD}
0 @>>> \hh(A,I) @>>> \Omega^1_{cl}(A) \otimes I @>{C^1 -
\Fr_{A^{[p]}}^* \otimes K_I}>> \Omega^1_{[p]}(A) \otimes I @>>> 0
\end{CD}
\end{equation}
of algebraic groups over $k$ (here we consider $k$-vector spaces
$\Omega^1_{cl}(A) \otimes I$ and $\Omega^1_{[p]}(A)$ as additive
algebraic groups, that is, we equip them with the scheme structure
of an affine space -- so that both become sums of several copies of
$\ga$). Note that by definition the group $\HH$ acts naturally on
the middle and the right-hand side term of this exact sequence,
hence also on the group $\hh(A,I)$.

\begin{remark}
As in Remark~\ref{not.lin}, we note that the map $\Fr_{A^{[p]}}^*
\otimes K_I$ is not $k$-linear, but only $k$-semilinear. Therefore
the commutative algebraic group $\hh(A,I)$ is {\em not} isomorphic
to the sum of several copies of $\ga$.  In fact, even the group
$\hh(A,I)(k)$ of its $k$-valued points has no natural structure of a
$k$-vector space. However, the algebraic group $\hh(A,I)$ is a
regular scheme over $k$.
\end{remark}

\begin{lemma}\label{ker.der}
Assume given a quantization base $B_0$ and a small extension $I
\subset B \to B_0$. Then the kernel of the surjective map $\G^B \to
\G^{B_0}$ is isomorphic to the algebraic group $\hh(A,I)$, and the
conjugation action of $\G^{B_0}$ on $\hh(A,I)$ is induced from the
natural action of $\HH$ on $\hh(A,I)$ via the reduction map
$\G^{B_0} \to \HH$.
\end{lemma}

\proof{} As in the proof of Proposition~\ref{small.ex}, the general
case immediately reduces to the case $B_0=k$, $\m_B=I$, $\m_B^2=0$
and $h=0$, so that $A_B$ is a restricted Poisson algebra over
$B$. Simplify notation by setting $K=K_I$. For any $k$-algebra $B'$,
the group $\hh(A,I)(B')$ is the subgroup in the group of $B'$-linear
Poisson derivations $D:A \otimes B' \to (A \otimes I) \otimes B'$
spanned by those derivations for which $\id \oplus D:A_B \otimes B'
\to A_B \otimes B'$ is compatible with restricted structures. This
group of derivations coincides with $\Omega^1_{cl}(A)(B')$ (unlike
$\hh(A,I)(B')$, it is not only a group, but also a $B'$-module). We
note that for any Poisson derivation $D:A \otimes B' \to A \otimes I
\otimes B'$, the composition $K \circ D$ is a Frobenius
derivation. By Theorem~\ref{obstr}, a map is compatible with a
restricted structure if and only if it preserves the associated
isomorphism
$$
C^1_*[\M_\Omega] \cong \Omega \cntrct [\theta_K].
$$
We leave it to the reader to check that in terms of $D$, this can be
rewritten as
$$
C^1(\Omega \cntrct D) = \Omega \cntrct (K \circ D),
$$
which coincides with \eqref{seq.3}.
\endproof

\begin{lemma}\label{splits}
Assume given a quantization base $B$ with maximal ideal $\m_B$ such
that $\m^2_B=0$. Then the natural surjective group homomorphism
$\G^B \to \HH$ obtained by reduction $\mod\m_B$ admits a splitting
$\HH \to \G^B$.
\end{lemma}

\proof{} If $h=0$, so that $A_B$ is a Poisson algebra, the claim is
obvious: $\HH$ naturally acts on $A_B = A \otimes_k B$. Denote by
$B'$ the associated graded quotient of the algebra $B$ with respect
to the $h$-adic filtration on $B$. We claim that $\G^B \cong
\G^{B'}$.

Indeed, since $k\langle 1 \rangle$ is injective in the category of
small Dieudonn\'e modules, there exists a splitting $I=\m_B \cong k
\cdot h \oplus (I/hI)$. Therefore there exists a restricted algebra
isomorphism $B \cong B'$. Let $\iota:B \to B$ be the involution that
fixes $B/hB$ and sends $h$ to $-h$, and let $\iota^*A^{opp}_B$ be
the {\em opposite} algebra to the $B$-quantization $A_B$ with the
$B$-module structure twisted by $\iota$, and the same Poisson
bracket and $p$-th power operation. Then it is easy to check that
$\iota^*A^{opp}_B$ is also a $B$-quantization of the algebra
$A$. Moreover, since $\iota$ is identical on $B_0 = B/ k \cdot h$,
both $A_B$ and $\iota^*A^{opp}_B$ are objects in $Q(A_0,B)$, where
$A_0 = A_B/A \cdot h$. By Proposition~\ref{q.uni}, there exists an
isomorphism $\iota:\iota^*A^{opp}_B \cong A_B$ identical on the
quotient $A_0$. In other words, we have a natural lifting $\iota:A_B
\to A_B$ of the map $\iota:B \to B$ to an anti-multiplicative map of
restricted quantized algebras. The map $\iota^2:A_B \to A_B$ is then
an automorphism of the algebra $A_B$ which is also identical on
$A_0$. Since $\iota^2$ commutes with $\iota$, while $\iota$ is
identical on $A_0$ and equal to $(-\id)$ on $A \cdot h \subset A_B$,
we must have $\iota^2=\id$, so that $\iota:A_B \to A_B$ is an
anti-involution. Taking its eigenspace decomposition, we obtain a
vector space isomorphism $A_B \cong A'_B$. It is easy to check that
this isomorphism is compatible with the Poisson bracket and the
$p$-th power operation, while the products in $A_B$ and $A'_B$ are
related by
$$
a \cdot b = ab + h\{a,b\}
$$
(here $a \cdot b$ is the product in $A_B$, and $ab$ is the product
in $A'_B$). We conclude that every $B$-linear automorphism of the
restricted quantized algebra $A_B$ fixes the product in $A'_B \cong
A_B$, and vice versa. Therefore we indeed have $\G_B \cong \G_B'$.
\endproof

\begin{remark}
In the case when $B=k[[h]]$ and $A_B = D$, the reduced Weyl algebra,
it is easy to compute the successive quotients in the natural
$h$-adic filtration on $\G=\G^B$ directly, without recourse to
restricted structures and deformation theory. Namely, one first
proves that the $D(h^{-1})$ is an Azumaya algebra over
$k((h))$. Then by the Scolem-Noether Theorem all $h$-linear
automorphisms of the algebra $D(h^{-1})$ are inner. It is an easy
exercise to check that an automorphism preserves $D \subset
D(h^{-1})$ and acts trivially on $A=D/hD$ if and only if it is given
by conjugation with an invertible element in $D \subset
D(h^{-1})$. One concludes that the kernel $\G^{> 0}$ of the natural
map $\G \to \HH$ is isomorphic to $D^*/k[[h]]^*$, where $D^*$ is the
group of invertible elements in $D$; the isomorphism in
Lemma~\ref{ker.der} immediately follows by applying the $d\log$
map. The problem with this approach is that it is not clear how to
generalize it to non-trivial quantization bases $B$, and it is
probably not possible at all to use it to prove Lemma~\ref{splits}.
\end{remark}

\section{Formal geometry.}

In order to obtain global results about quantizations, we use the
technique of the so-called {\em formal geometry}. This is similar to
the approach in \cite[3.1]{BK} in the characteristic $0$ case;
however, the situation in positive characteristic is somewhat
simpler, and in particular, all the algebraic groups we use are of
finite type.

\subsection{The setup.}

Fix a finite-dimensional vector space $W$ over $k$, and let $A =
k[W^*]/\m^{[p]}$ be the commutative $k$-algebra associated to $W$ in
Subsection~\ref{drb.sub}. Recall that we have the algebraic group
$\Aut A$ of all automorphisms of the commutative algebra $A$ and the
reduced subgroup $(\Aut A)_0 \subset \Aut A$ of automorphisms which
preserve the maximal ideal $\m \subset A$. The quotient $(\Aut
A)/(\Aut A)_0$ is naturally identified with the spectrum $\Spec A$
of the algebra $A$ -- indeed, the action of $\Aut A$ on $\Spec A$ is
transitive, and $(\Aut A)_0 \subset \Aut A$ is the stabilizer of the
closed point.

Formal geometry is based on the following observation. Let $X$ be a
manifold over $k$, -- that is, a quasiregular reduced scheme $X/k$
of finite type, -- and assume that $\rk(\T(X))=\dim W$. Then there
exists a canonical $(\Aut A)_0$-torsor $\M_X$ over $X$. Namely, for
any affine $S/k$ we define
$$
\M_X(S) = \Hom_{\text{\'et}}(\Spec A \times S,X \times S),
$$
the set of all unramified maps $\phi:\Spec A \times S \to X \times
S$ of schemes over $S$. To check that the functor $\M_X(S)$ is
represented by a faithfully flat scheme over $X$ is an easy exercise
left to the reader. The group $\Aut A$ acts on $\M_X$ through its
action on $\Spec A$; a second easy check shows that the action is
effective, and that $X \cong \M_X/(\Aut A)_0$.

Informally, $\M_X$ is the scheme of all pairs $\langle x, \phi
\rangle$ of a point $x \in X$ and a surjective map $\calo_{X,x} \to
A$ -- a ``coordinate system in the Frobenius neighborhood of $x \in
X$''. We will call $\M_X$ the {\em torsor of Frobenius frames} on
$X$. The projection $\M_X \to X = \M_X/(\Aut A)_0$ corresponds to
forgetting the coordinate system. We note that the torsor of
Frobenius frames is locally trivial in Zariski topology on $X$.

The action of the whole group $\Aut A$ on $X$ is also effective, and
the quotient $\M_X/(\Aut A)$ is naturally identified with $X^{[p]}$,
that is, $X$ with structure sheaf $\calo_X^p \subset \calo_X$. All
in all, we have a diagram
$$
\begin{CD}
\M_X @>>> X @>{\Fr}>> X^{[p]}.
\end{CD}
$$
Here the left-hand side map is the quotient map $\M_X \to \M_X/(\Aut
A)_0 \cong X$, the composition is the quotient map $\M_X \to
\M_X/(\Aut A)_0 \cong X^{[p]}$, and the right-hand side map is the
Frobenius map of the manifold $X$. The scheme $\M_X$ can be treated
either as an $(\Aut A)_0$-torsor over $X$, or as an $\Aut A$-torsor
over $X^{[p]}$.

Having said this, we can now describe various differential-geometric
structures on $X$ in terms of reductions of the torsor $\M_X$ to
various subgroups in $\Aut A$ -- or, more generally, its liftings to
various algebraic groups $G$ equipped with a map $G \to \Aut A$.

\begin{defn}\label{G.st.defn}
For any algebraic group $G$ over $k$ equipped with a map $\rho:G \to \Aut
A$, a {\em $G$-structure} on the manifold $X$ is a $G$-torsor $\M_G$
over $X^{[p]}$ equipped with a map $\M_G \to \M_X$ compatible with
$\rho$.
\end{defn}

For any $k$-vector space $V$ equipped with an algebraic action of the
group $\Aut A$, we have the associated bundle $\loc(\M_X,V)$; this
defines an exact {\em localization functor}
$$
\loc:(\Aut A)\fmod \to \Coh(X^{[p]})
$$
from the category of finite-dimensional representations of the group
$\Aut A$ to the category of coherent sheaves of $\calo_X^p$-modules
on $X^{[p]}$. Essentially this is just the flat descent functor for
the $(\Aut A)$-torsor $\M_X$. For example, the localization
$\loc(k)$ of the trivial $(\Aut A)$-module $k$ is isomorphic to the
structure sheaf $\calo_X^p$. More generally, assume given a
commutative algebraic group $H$ over $k$ -- equivalently, a sheaf of
abelian groups on $\Spec k$ in the flat topology -- and assume that
the group $\Aut A$ acts on $H$ in an algebraic way. Then flat
descent allows to define a sheaf $\loc(\M_G,H)$ of abelian groups on
$X^{[p]}$ in flat topology.  This works just as well for any
$G$-structure $\M_G$ on $X$ -- we have a localization functor
$\loc(\M_G,-)$ from the category of abelian algebraic groups over
$k$ equipped with a $G$-action to the category of sheaves of abelian
groups on $X^{[p]}$ in the flat topology. We recall that both
these categories are abelian; the localization functor is exact.

Using localization, one defines obstruction theory for
$G$-structures in the following way. Assume given an algebraic group
$G$ over $k$ equipped with a map $G \to \Aut A$, and another
algebraic group $G_1$ over $k$ equipped with a surjective map $G_1
\to G$. Moreover, assume that the kernel $H$ of the map $G_1 \to G$
is an abelian algebraic group. Then $G$ acts by conjugation on the
group $H$, and this action is algebraic. The extension $G_1 \to G$
is completely defined by a cohomology class $c \in H^2(G,H)$ (in
particular, $G_1$ coincides with the semi-direct product $H \rtimes
G$ if and only if $c=0$). Assume given a $G$-structure $\M$ on a
manifold $X$ of dimension $\dim X = \dim W$. Then we can apply
localization to the class $c \in H^2(G,H) = \Ext^2(\ga,H)$ and
obtain a cohomology class
$$
\loc(\M,c) \in \Ext^2(\calo_X^p,\loc(\M,V)) =
H^2(X^{[p]},\loc(\M,H)).
$$
The following is the standard obstruction theory statement; for the
proof see e.g. \cite{gir}.

\begin{lemma}\label{brau}
The $G$-torsor $\M$ admits a restriction to a $G_1$-torsor $\M_1$ if
and only if $\loc(\M,c)=0$. Moreover, if this happens, then for any
two such restrictions $\M_1$, $\M_1'$ we have a well-defined class
$[\M_1 - \M_1']) \in H^1(X^{[p]},\loc(\M,H))$, and $\M_1 \cong
\M_1'$ if and only if $[\M_1 - \M_1']) = 0$.\endproof
\end{lemma}

\subsection{The proofs.}\label{azu.sub}\label{step.pf.sub}
We can now prove Proposition~\ref{step}, Proposition~\ref{level.1}
and Proposition~\ref{azu}.

\begin{lemma}\label{q=G}
Let $B$ be a quantization base over the field $k$, let $X$ be a
symplectic manifold of dimension $\dim X = \dim W$, and let $\G^B$
be the automorphism group of the regular $B$-quantization $A_B$ of
the algebra $A$ introduced in Proposition~\ref{q.uni}. Then there is
a natural equivalence between the category of regular
$B$-quantizations of $X$ and the category of $\G^B$-structures on
$X$. Moreover, every $\G^B$ structure $\M_G$ on $X$ is locally
trivial in \'etale topology.
\end{lemma}

\proof{} Given a $\G^B$-structure $\M_B$ on $X$, we obtain a regular
$B$-quantization $\calo_B$ as
$$
\calo_B = \loc(\M_B,A_B).
$$
Conversely, given a regular $B$-quantization $\calo_B$, we take
$\M_B$ to be the set of pairs,
$$
\M = \left\{ \langle x, \phi \rangle \mid x \in
X^{[p]},\phi:\left(\calo_B\right)_x \to A_B \right\},
$$
where $x \in X^{[p]}$ is a point, and $\phi$ is an isomorphism
between the stalk $\left(\calo_B\right)_x$ of the $\calo_X^p$-module
sheaf $\calo_B$ at the point $x$ and the standard restricted
quantized algebra $A_B$. One checks easily that $\M_B$ has a natural
algebraic structure; the group $\G^B$ acts on $\M_B$ by $g(\langle
x,\phi\rangle) = \langle x, g \circ \phi \rangle$, this action is
effective, and the natural projection $\M_B \to X^{[p]}$ identifies
$\M_B/\G^B$ with an open subset in $X^{[p]}$. Finally, by
Proposition~\ref{q.uni} the algebra $\left(\calo_B\right)_x$ is
isomorphic to $A_B$ for any $x \in X^{[p]}$, so that $\M_B/\G^B$ is
the whole $X^{[p]}$, and $\M_B$ is indeed a $\G^B$-structure on
$X$. To prove the second claim, we note that by
Proposition~\ref{ker.der} the kernel of the map $\G^B \to \Aut A$ is
an iterated extension of smooth commutative algebraic groups of the
form $\hh\langle A, I \rangle$. Therefore the torsor $\M_B$ over
$\M_X$ is locally trivial in \'etale topology (see e.g. \cite[III,
\S 3, Theorem 3.9]{mln}). But $\M_X$ is in turn locally trivial in
\'etale -- and even in Zariski -- topology on $X$.
\endproof

\proof[Proof of Proposition~\ref{step} and
Proposition~\ref{level.1}.]\label{step.pf} Proposition~\ref{step} is
a direct combination of Lemma~\ref{q=G}, Lemma~\ref{brau} and
Proposition~\ref{ker.der}. To deduce Proposition~\ref{level.1}, it
suffices to add Lemma~\ref{splits} to the mix.
\endproof

\proof[Proof of Proposition~\ref{azu}.]\label{azu.pf} The first
claim follows from Lemma~\ref{loc.azu}. To prove the second claim,
we note that by Lemma~\ref{q=G}, the algebra sheaf $\calo_h$ is the
associated bundle to a $\G$-torsor $\M_h$ on $X^\tw$, where
$\G=\Aut(D)$ is the automorphism group of the reduced Weyl algebra
$D$. The group $\G$ is in fact obtained by Weil restriction of
scalars from a group $\wh{\G}$ over the Taylor power series algebra
$k[[h]]$, and the torsor $\M_h$ is obtained by restriction of
scalars from a torsor $\wh{\M}_h$ over the scheme $\wh{X}$. Over
$\overline{X} = \wh{X} \setminus X^\tw \subset \wh{X}$, the torsor
$\wh{\M}_h$ is a torsor over the group $\overline{\G}=\wh{\G}
\otimes_{k[[h]]} k((h))$. However, we have
$$
\overline{\G} \cong \Aut(D\otimes_{k[[h]]} k((h))) \cong
PGL\left(p^{\frac{\dim X}{2}},k((h))\right),
$$
and the torsor $\wh{\M}_h$ over $\overline{X}$ becomes exactly the
principal bundle associated to the Azumaya algebra
$\calo_h(h^{-1})$. Its class $\overline{\Br([\calo_h])} \in
\Br(\overline{X})$ is represented by a gerb bound by
$\calo_{\overline{X}}^*$ which is associated to the torsor
$\wh{\M}_h$ and the standard central extension
\begin{equation}\label{ccn}
\begin{CD}
k((h))^* @>>> GL\left(p^{\frac{\dim X}{2}},k((h))\right) @>>>
PGL\left(p^{\frac{\dim X}{2}},k((h))\right).
\end{CD}
\end{equation}
Since after restriction to $\wh{\G} \subset \overline{\G}$ this
extension reduces to a central extension of the group $\wh{\G}$ by
$k[[h]]^*$, the gerb representing $\overline{\Br([\calo_h])}$
extends to a gerb on $\wh{X}$ bound by $\calo_{\wh{X}}^*$, so that
the class $\overline{\Br([\calo_h])}$ comes from a class
$\wh{\Br([\calo_h])} \in \Br(\wh{X})$. However, by assumption
$H^2(X^\tw,\calo_X^p)=0$; considering the spectral sequence
associated to the $h$-adic filtration on $\calo_{\wh{X}}^*$, we
conclude that the restriction map
$$
\Br(\wh{X}) \to \Br(X^\tw)
$$
is an isomorphism. To finish the proof, it suffices to check that
the restriction of the class $\wh{\Br([\calo_h])}$ to $X^\tw \subset
\wh{X}$ indeed coincides with $\Br([\calo_h])$. This immediately
follows from the commutative diagram
$$
\begin{CD}
k[[h]]^* @>>>      \wt{\G}       @>>> \G^{k[[h]]}\\
@VVV @VVV @VVV\\
k^*      @>>>  A^*\rtimes\G_0    @>>> \G^{k[h]/h^2}
\end{CD},
$$ 
where the rows are central extensions of algebraic groups, and
$\wt{\G}$ is the extension induced by \eqref{ccn}.
\endproof

\begin{remark}
As the above proof shows, there is no canonical way to extend the
Azumaya algebra $\calo_h(h^{-1})$ itself to an Azumaya algebra on
$\wh{X}$ -- it is only its class in the Brauer group that has such
an extension. The reason for that is clear from the group theoretic
description: $\Aut(D)$ and $PGL(p^{\frac{\dim X}{2}},k[[h]])$ are
two different parabolic subgroups in the loop group
$PGL(p^{\frac{\dim X}{2}},k((h)))$, and both are maximal
parabolic. Extending the Azumaya algebra would correspond to finding
an embedding $\Aut(D) \subset PGL(p^{\frac{\dim X}{2}},k[[h]]))$;
since both are maximal parabolic, such an embedding does not exist.
\end{remark}

\bigskip

\noindent
{\sc Massachusets Institute of Technology\\
Cambridge, MA, USA\\
\medskip
and\\
\smallskip
Steklov Math Institute\\
Moscow, USSR}

\bigskip

\noindent
{\em E-mail addresses\/}: {\tt bezrukav@math.mit.edu}\\
\phantom{{\em E-mail addresses\/}: }{\tt kaledin@mccme.ru}

\end{document}